\documentclass[12pt, reqno]{amsart}
\usepackage{txfonts}      %{article} was 12pt latex e
\usepackage{amssymb}
\usepackage{eucal}
\usepackage{amsmath}
\usepackage{amsthm}
\usepackage{amscd}
\usepackage[dvips]{color}
\usepackage{multicol}
\usepackage[all]{xy}           %xypic macro for latex2.09
\usepackage{color}
\usepackage{colordvi}
\usepackage{xspace}
\usepackage{tikz}
\usepackage{ifpdf}
\ifpdf
\usepackage[colorlinks,final,backref=page,hyperindex]{hyperref}
\else
\usepackage[colorlinks,final,backref=page,hyperindex,hypertex]{hyperref}
\fi

\newcommand {\emptycomment}[1]{} %to remove paragraphs

\newcommand{\nc}{\newcommand}
\newcommand{\delete}[1]{}

	\nc{\mlabel}[1]{\label{#1}}  % Use this to suppress names
	\nc{\mcite}[1]{\cite{#1}}  % Use this to suppress names
	\nc{\mref}[1]{\ref{#1}}  % Use this to suppress names
	\nc{\meqref}[1]{\eqref{#1}} % Use this to suppress names
	\nc{\mbibitem}[1]{\bibitem{#1}} % Use this to show number
\delete{
	\nc{\mlabel}[1]{\label{#1}  % Use the next two lines to show names
		{\hfill \hspace{1cm}{\bf{{\ }\hfill(#1)}}}}
	\nc{\mcite}[1]{\cite{#1}{{\bf{{\ }(#1)}}}}  % Use this lines to show names
	\nc{\mref}[1]{\ref{#1}{{\bf{{\ }(#1)}}}}  % Use this lines to show names
	\nc{\meqref}[1]{\eqref{#1}{{\bf{{\ }(#1)}}}} % Use this lines to show names
	\nc{\mbibitem}[1]{\bibitem[\bf #1]{#1}} % Use this to show name
}

\newtheorem{thm}{Theorem}[section]
\newtheorem{lem}[thm]{Lemma}
\newtheorem{cor}[thm]{Corollary}
\newtheorem{pro}[thm]{Proposition}
\theoremstyle{definition}
\newtheorem{defi}[thm]{Definition}
\newtheorem{ex}[thm]{Example}
\newtheorem{rmk}[thm]{Remark}

\nc{\tred}[1]{\textcolor{red}{#1}}
\nc{\tblue}[1]{\textcolor{blue}{#1}}
\nc{\tgreen}[1]{\textcolor{green}{#1}}
\nc{\tpurple}[1]{\textcolor{purple}{#1}}
\nc{\btred}[1]{\textcolor{red}{\bf #1}}
\nc{\btblue}[1]{\textcolor{blue}{\bf #1}}
\nc{\btgreen}[1]{\textcolor{green}{\bf #1}}
\nc{\btpurple}[1]{\textcolor{purple}{\bf #1}}

\nc{\rk}{\mathrm{r}}

\nc{\tforall}{\text{ for all }}

\nc{\svec}[2]{{\tiny\left(\begin{matrix}#1\\
			#2\end{matrix}\right)\,}}  % column vector
\nc{\ssvec}[2]{{\tiny\left(\begin{matrix}#1\\
			#2\end{matrix}\right)\,}} % subscript column vector
\nc{\typeI}{local cocycle $3$-Lie bialgebra\xspace}
\nc{\typeIs}{local cocycle $3$-Lie bialgebras\xspace}
\nc{\typeII}{double construction $3$-Lie bialgebra\xspace}
\nc{\typeIIs}{double construction $3$-Lie bialgebras\xspace}
\nc{\bia}{{$\mathcal{P}$-bimodule ${\bf k}$-algebra}\xspace}
\nc{\bias}{{$\mathcal{P}$-bimodule ${\bf k}$-algebras}\xspace}
\nc{\rmi}{{\mathrm{I}}}
\nc{\rmii}{{\mathrm{II}}}
\nc{\rmiii}{{\mathrm{III}}}
\nc{\OT}{constant $\theta$-}
\nc{\T}{$\theta$-}
\nc{\IT}{inverse $\theta$-}

\nc{\asi}{ASI\xspace}
\nc{\adm}{admissible\xspace}
\nc{\qadm}{$Q$-admissible\xspace}
\nc{\aybe}{AYBE\xspace}
\nc{\aybes}{AYBEs\xspace}
\nc{\admset}{\{\pm x\}\cup (-x+K^{\times}) \cup K^{\times} x^{-1}}
\nc{\dualrep}{gives a dual representation\xspace}
\nc{\admt}{admissible to\xspace}
\nc{\opa}{\cdot_A}
\nc{\opb}{\cdot_B}

\nc{\dfl}{\succ} \nc{\dfr}{\prec}
\nc{\disp}[1]{\displaystyle{#1}}
\nc{\bin}[2]{ (_{\stackrel{\scs{#1}}{\scs{#2}}})}  %binomial coeff
\nc{\binc}[2]{ \left (\!\! \begin{array}{c} \scs{#1}\\
		\scs{#2} \end{array}\!\! \right )}  %binomial coeff
\nc{\bincc}[2]{  \left ( {\scs{#1} \atop
		\vspace{-.5cm}\scs{#2}} \right )}  %binomial coeff
\nc{\sarray}[2]{\begin{array}{c}#1 \vspace{.1cm}\\ \hline
		\vspace{-.35cm} \\ #2 \end{array}}
\nc{\bs}{\bar{S}} \nc{\dcup}{\stackrel{\bullet}{\cup}}
\nc{\dbigcup}{\stackrel{\bullet}{\bigcup}} \nc{\etree}{\big |}
\nc{\la}{\longrightarrow} \nc{\fe}{\'{e}} \nc{\rar}{\rightarrow}
\nc{\dar}{\downarrow}
\nc{\ot}{\otimes} \nc{\sot}{{\scriptstyle{\ot}}}
\nc{\otm}{\overline{\ot}}
\nc{\ora}[1]{\stackrel{#1}{\rar}}
\nc{\ola}[1]{\stackrel{#1}{\la}}%${\Bbb Z}$
\nc{\pltree}{\calt^\pl}
\nc{\epltree}{\calt^{\pl,\NC}}
\nc{\rbpltree}{\calt^r}
\nc{\scs}[1]{\scriptstyle{#1}} \nc{\mrm}[1]{{\rm #1}}
\nc{\dirlim}{\displaystyle{\lim_{\longrightarrow}}\,}
\nc{\invlim}{\displaystyle{\lim_{\longleftarrow}}\,}
\nc{\sha}{{\mbox{\cyr X}}}  %used to be \cyr
\nc{\ncsha}{{\mbox{\cyr X}^{\mathrm NC}}} \nc{\ncshao}{{\mbox{\cyr
			X}^{\mathrm NC,\,0}}}
\nc{\shpr}{\diamond}    %Shuffle product
\nc{\shprm}{\overline{\diamond}}    %Shuffle product
\nc{\shpro}{\diamond^0}    %Shuffle product
\nc{\shprr}{\diamond^r}     %product on controlled trees
\nc{\shpra}{\overline{\diamond}^r}
\nc{\shpru}{\check{\diamond}} \nc{\catpr}{\diamond_l}
\nc{\rcatpr}{\diamond_r} \nc{\lapr}{\diamond_a}
\nc{\sqcupm}{\ot}
\nc{\lepr}{\diamond_e} \nc{\vep}{\varedhlon} \nc{\labs}{\mid\!}
\nc{\rabs}{\!\mid} \nc{\hsha}{\widehat{\sha}}
\nc{\lsha}{\stackrel{\leftarrow}{\sha}}
\nc{\rsha}{\stackrel{\rightarrow}{\sha}} \nc{\lc}{\lfloor}
\nc{\rc}{\rfloor}
\nc{\tpr}{\sqcup}
\nc{\nctpr}{\vee}
\nc{\plpr}{\ob}
\nc{\rbplpr}{\bar{\plpr}}
\nc{\sqmon}[1]{\langle #1\rangle}
\nc{\forest}{\calf}
\nc{\altx}{\Lambda_X} \nc{\vecT}{\vec{T}} \nc{\onetree}{\bullet}
\nc{\Ao}{\check{A}}
\nc{\seta}{\underline{\Ao}}
\nc{\deltaa}{\overline{\delta}}
\nc{\trho}{\tilde{\rho}}

\nc{\rpr}{\rhd}
\nc{\dpr}{{\tiny\diamond}}
\nc{\rprpm}{{\rpr}}

\nc{\mmbox}[1]{\mbox{\ #1\ }} \nc{\ann}{\mrm{ann}}
\nc{\Aut}{\mrm{Aut}} \nc{\can}{\mrm{can}}
\nc{\twoalg}{{two-sided algebra}\xspace}
\nc{\colim}{\mrm{colim}}
\nc{\Cont}{\mrm{Cont}} \nc{\rchar}{\mrm{char}}
\nc{\cok}{\mrm{coker}} \nc{\dtf}{{R-{\rm tf}}} \nc{\dtor}{{R-{\rm
			tor}}}

\nc{\depth}{{\mrm d}}
\nc{\Div}{{\mrm Div}} \nc{\End}{\mrm{End}} \nc{\Ext}{\mrm{Ext}}
\nc{\Fil}{\mrm{Fil}} \nc{\Frob}{\mrm{Frob}} \nc{\Gal}{\mrm{Gal}}
\nc{\GL}{\mrm{GL}} \nc{\Hom}{\mrm{Hom}} \nc{\hsr}{\mrm{H}}
\nc{\hpol}{\mrm{HP}} \nc{\id}{\mrm{id}} \nc{\im}{\mrm{im}}
\nc{\incl}{\mrm{incl}} \nc{\length}{\mrm{length}}
\nc{\LR}{\mrm{LR}} \nc{\mchar}{\rm char} \nc{\NC}{\mrm{NC}}
\nc{\mpart}{\mrm{part}} \nc{\pl}{\mrm{PL}}
\nc{\ql}{{\QQ_\ell}} \nc{\qp}{{\QQ_p}}
\nc{\rank}{\mrm{rank}} \nc{\rba}{\rm{RBA }} \nc{\rbas}{\rm{RBAs }}
\nc{\rbpl}{\mrm{RBPL}}
\nc{\rbw}{\rm{RBW }} \nc{\rbws}{\rm{RBWs }} \nc{\rcot}{\mrm{cot}}
\nc{\rest}{\rm{controlled}\xspace}
\nc{\rdef}{\mrm{def}} \nc{\rdiv}{{\rm div}} \nc{\rtf}{{\rm tf}}
\nc{\rtor}{{\rm tor}} \nc{\res}{\mrm{res}} \nc{\SL}{\mrm{SL}}
\nc{\Spec}{\mrm{Spec}} \nc{\tor}{\mrm{tor}} \nc{\Tr}{\mrm{Tr}}
\nc{\mtr}{\mrm{sk}}

\nc{\ab}{\mathbf{Ab}} \nc{\Alg}{\mathbf{Alg}}
\nc{\Algo}{\mathbf{Alg}^0} \nc{\Bax}{\mathbf{Bax}}
\nc{\Baxo}{\mathbf{Bax}^0} \nc{\RB}{\mathbf{RB}}
\nc{\RBo}{\mathbf{RB}^0} \nc{\BRB}{\mathbf{RB}}
\nc{\Dend}{\mathbf{DD}} \nc{\bfk}{{\bf k}} \nc{\bfone}{{\bf 1}}
\nc{\base}[1]{{a_{#1}}} \nc{\detail}{\marginpar{\bf More detail}
	\noindent{\bf Need more detail!}
	\svp}
\nc{\Diff}{\mathbf{Diff}} \nc{\gap}{\marginpar{\bf
		Incomplete}\noindent{\bf Incomplete!!}
	\svp}
\nc{\FMod}{\mathbf{FMod}} \nc{\mset}{\mathbf{MSet}}
\nc{\rb}{\mathrm{RB}} \nc{\Int}{\mathbf{Int}}
\nc{\Mon}{\mathbf{Mon}}
\nc{\remarks}{\noindent{\bf Remarks: }}
\nc{\OS}{\mathbf{OS}} %free operated semigroup
\nc{\Rep}{\mathbf{Rep}}
\nc{\Rings}{\mathbf{Rings}} \nc{\Sets}{\mathbf{Sets}}
\nc{\DT}{\mathbf{DT}}

\nc{\BA}{{\mathbb A}} \nc{\CC}{{\mathbb C}} \nc{\DD}{{\mathbb D}}
\nc{\EE}{{\mathbb E}} \nc{\FF}{{\mathbb F}} \nc{\GG}{{\mathbb G}}
\nc{\HH}{{\mathbb H}} \nc{\LL}{{\mathbb L}} \nc{\NN}{{\mathbb N}}
\nc{\QQ}{{\mathbb Q}} \nc{\RR}{{\mathbb R}} \nc{\BS}{{\mathbb{S}}} \nc{\TT}{{\mathbb T}}
\nc{\VV}{{\mathbb V}} \nc{\ZZ}{{\mathbb Z}}

\nc{\calao}{{\mathcal A}} \nc{\cala}{{\mathcal A}}
\nc{\calc}{{\mathcal C}} \nc{\cald}{{\mathcal D}}
\nc{\cale}{{\mathcal E}} \nc{\calf}{{\mathcal F}}
\nc{\calfr}{{{\mathcal F}^{\,r}}} \nc{\calfo}{{\mathcal F}^0}
\nc{\calfro}{{\mathcal F}^{\,r,0}} \nc{\oF}{\overline{F}}
\nc{\calg}{{\mathcal G}} \nc{\calh}{{\mathcal H}}
\nc{\cali}{{\mathcal I}} \nc{\calj}{{\mathcal J}}
\nc{\call}{{\mathcal L}} \nc{\calm}{{\mathcal M}}
\nc{\caln}{{\mathcal N}} \nc{\calo}{{\mathcal O}}
\nc{\calp}{{\mathcal P}} \nc{\calq}{{\mathcal Q}} \nc{\calr}{{\mathcal R}}
\nc{\calt}{{\mathcal T}} \nc{\caltr}{{\mathcal T}^{\,r}}
\nc{\calu}{{\mathcal U}} \nc{\calv}{{\mathcal V}}
\nc{\calw}{{\mathcal W}} \nc{\calx}{{\mathcal X}}
\nc{\CA}{\mathcal{A}}

\nc{\fraka}{{\mathfrak a}} \nc{\frakB}{{\mathfrak B}}
\nc{\frakb}{{\mathfrak b}} \nc{\frakd}{{\mathfrak d}}
\nc{\oD}{\overline{D}}
\nc{\frakF}{{\mathfrak F}} \nc{\frakg}{{\mathfrak g}}
\nc{\frakm}{{\mathfrak m}} \nc{\frakM}{{\mathfrak M}}
\nc{\frakMo}{{\mathfrak M}^0} \nc{\frakp}{{\mathfrak p}}
\nc{\frakS}{{\mathfrak S}} \nc{\frakSo}{{\mathfrak S}^0}
\nc{\fraks}{{\mathfrak s}} \nc{\os}{\overline{\fraks}}
\nc{\frakT}{{\mathfrak T}} \nc{\frakO}{{\mathfrak O}}
\nc{\oT}{\overline{T}}
\nc{\frakX}{{\mathfrak X}} \nc{\frakXo}{{\mathfrak X}^0}
\nc{\frakx}{{\mathbf x}}
\nc{\frakTx}{\frakT}      %All rooted trees, correspond to \ncsha(X)
\nc{\frakTa}{\frakT^a}        % rooted trees for \ncsha(A)
\nc{\frakTxo}{\frakTx^0}   % rooted trees for \ncshao(X)
\nc{\caltao}{\calt^{a,0}}   % rooted trees for \ncshao(A)
\nc{\ox}{\overline{\frakx}} \nc{\fraky}{{\mathfrak y}}
\nc{\frakz}{{\mathfrak z}} \nc{\oX}{\overline{X}}

\font\cyr=wncyr10

\nc{\al}{\alpha}
\nc{\lam}{\lambda}
\nc{\lr}{\longrightarrow}

\nc{\pr}{\prec}
\nc{\su}{\succ}
\nc{\prp}{\prec_P}
\nc{\scp}{\succ_P}
\nc{\pra}{\prec_A}
\nc{\sca}{\succ_A}
\nc{\prb}{\prec_B}
\nc{\scb}{\succ_B}
\nc{\nea}{\nearrow} % WARNING: \ne is "not equal to" sign
\nc{\se}{\searrow}
\nc{\sw}{\swarrow}
\nc{\nw}{\nwarrow}
\nc{\nep}{\ne^{op}}
\nc{\sep}{\se^{op}}
\nc{\swp}{\sw^{op}}
\nc{\nwp}{\nw^{op}}
\nc{\east}{\succ}
\nc{\west}{\prec}
\nc{\north}{\wedge}
\nc{\south}{\vee}
\nc{\eastt}{\east^t}
\nc{\westt}{\west^t}
\nc{\northt}{\north^t}
\nc{\ob}{{\ \begin{picture}(-1,1)(-1,-3)\circle*{3}\end{picture}\ \,}}%{\cdot}

\nc{\dftimes}{\blacksquare}

\nc{\asib}{antisymmetric infinitesimal bialgebra\xspace}
\nc{\ASIb}{ASI bialgebra\xspace}
\nc{\asibs}{antisymmetric infinitesimal bialgebras\xspace}
\nc{\ASIbs}{ASI bialgebras\xspace}

\nc{\lsa}{\ell_{\prec_A}}
\nc{\rsa}{r_{\prec_A}}
\nc{\lusa}{\ell_{{\succ_A}}}
\nc{\rusa}{r_{{\succ_A}}}
\nc{\lsb}{\ell_{\prec_B}}
\nc{\rsb}{r_{\prec_B}}
\nc{\lusb}{\ell_{{\succ_B}}}
\nc{\rusb}{r_{{\succ_B}}}
\nc{\at}{\cdot_d}

\nc{\wt}{of weight $\lambda$\xspace}
\nc{\pa}{\partial}

\nc{\NYBE}{NS-YBE\xspace}

\nc{\frakr}{{\mathfrak r}}
\nc{\vsa}{\vspace{-.1cm}}
\nc{\vsb}{\vspace{-.2cm}}
\nc{\vsc}{\vspace{-.3cm}}
\nc{\vsd}{\vspace{-.4cm}}
\nc{\vse}{\vspace{-.5cm}}

\nc{\sh}[1]{\textcolor{purple}{Shanghua:#1}}
\nc{\yz}[1]{\textcolor{green}{Yizhen:#1}}
\nc{\li}[1]{\textcolor{red}{#1}}
\nc{\lir}[1]{\textcolor{blue}{Li:#1}}
\begin{document}
	
	\title[NS-bialgebras]{
NS-bialgebras, the NS-Yang-Baxter equation and dendriform D-bialgebras}

\author{Shanghua Zheng}
\address{School of Mathematics and Statistics, Jiangxi Normal University, Nanchang, Jiangxi 330022, China}
\email{zhengsh@jxnu.edu.cn}

\author{Chenyu Qiu}
\address{School of Mathematics and Statistics, Jiangxi Normal University, Nanchang, Jiangxi 330022, China}
\email{qiuchenyu@jxnu.edu.cn}

\author{Xigou Zhang}
\address{School of Mathematics and Statistics, Jiangxi Normal University, Nanchang, Jiangxi 330022, China}
\email{xyzhang@jxnu.edu.cn}

	%========================================================================
	%\date{\today}
	%========================================================================
\begin{abstract}
In this paper, we develop a bialgebra theory for NS-algebras via the Manin triple approach. We first study representations of NS-algebras and establish an equivalence among NS-bialgebras, matched pairs and Manin triples of NS-algebras. We then introduce the NS-Yang-Baxter equation (NS-YBE) and show that its symmetric solutions induce triangular NS-bialgebras. We further investigate the relationship between symmetric solutions of the NS-YBE and $\mathcal O$-operators of NS-algebras, and characterize pre-NS-algebras in terms of invertible $\mathcal O$-operators.
As an application, we investigate the relation between NS-bialgebras and dendriform D-bialgebras via double constructions of Connes cocycles. In particular, we show that every NS-bialgebra induces a dendriform D-bialgebra structure on the direct sum of two copies of its underlying vector space.
\end{abstract}

	\makeatletter
	\@namedef{subjclassname@2020}{\textup{2020} Mathematics Subject Classification}
	\makeatother
	\subjclass[2020]{
		16T10, %Bialgebras
		17B38, %Yang-Baxter equations and Rota-Baxter operators
		17A30, %Nonassociative algebras satisfying other identities
		%17D25, %Lie-admissible algebras
		16T15, %Coalgebras and comodules; corings
		17A36. %Automorphisms, derivations, other operators (nonassociative rings and algebras)
		%18M65, %Non-symmetric operads, multicategories, generalized multicategories
		%12H05, %Differential algebra
		%16S10, %Rings determined by universal properties (free algebras, coproducts, adjunction of inverses, etc.)
		%16W99, %Rings and algebras with additional structure/none of the above
		%16T30,  %Connections with combinatorics
		%57R56
			}
	
	\keywords{NS-algebra, bialgebra, Nijenhuis operator, Manin triple, 
		Yang-Baxter equation,
		$\mathcal{O}$-operator,
		dendriform D-bialgebra}
	\maketitle
	
	\vspace{-1cm}
	
	\tableofcontents
	
	\vspace{-1cm}
	
	\allowdisplaybreaks

\section{Introduction}

NS-algebras were introduced by Leroux~\cite{Ler04}
in connection with Nijenhuis operators on associative
algebras. They provide a splitting of associative
multiplication into three binary operations and are
closely related to dendriform algebras.
Lei and Guo~\cite{LG12} further investigated
the relationships among Nijenhuis algebras,
NS-algebras and N-dendriform algebras, including
constructions of free Nijenhuis algebras and
universal enveloping Nijenhuis algebras of
NS-algebras.

The theory of NS-algebras is also connected
with Rota-Baxter type operatorsand twisting
constructions~\cite{Uch10, ZGGS21}. Uchino~\cite{Uch08}
studied relationships among Rota-Baxter
operators, dendriform structures and NS-algebras.
In particular, twisted Rota-Baxter operators
provide another source of NS-algebra structures.
Das~\cite{Das22a} developed a cohomology theory
for NS-algebras and investigated their
deformations, while twisted Rota-Baxter
families and the associated NS-family algebras
were studied in~\cite{Das22b}.
More recently, generalized NS-algebras
were introduced in a broader algebraic
framework~\cite{OPV25}.
Nijenhuis deformations of Poisson and
$F$-manifold algebras, together with related
NS-type structures, were investigated
in~\cite{BD26}.

The NS-algebra framework has also been
extended to several nonassociative settings.
For Lie algebras, twisted Rota-Baxter
and Reynolds operators lead to NS-Lie
structures~\cite{Das21}.
Related constructions for Leibniz algebras,
pre-Lie algebras and higher Lie-type
algebras have been developed
in~\cite{DG21,GMM23,GZ23,HS21,LW23}.
Hom-associative analogues involving twisted
Rota-Baxter family operators and their
deformations were studied in~\cite{Teng25}.
These developments illustrate the connections
between NS-type structures, operator identities
and deformation theory.

Bialgebra structures provide a framework
for studying the compatibility between
algebraic operations and their dual
coalgebraic counterparts.
The theory of Lie bialgebras and the
classical Yang-Baxter equation, initiated
by Drinfeld~\cite{D1983}, has motivated
the development of bialgebra theories
associated with various nonassociative
algebras.
Related algebraic structures also arise
in the study of Hamiltonian
operators~\cite{GD79}.
For Leibniz algebras, Tang and
Sheng~\cite{TS22} established the
equivalence among Leibniz bialgebras,
matched pairs and Manin triples, and
investigated relative Rota-Baxter
operators and the classical Leibniz
Yang-Baxter equation.
A corresponding approach to Nijenhuis
Lie bialgebras, including cohomological
aspects, was developed in~\cite{Das25}.
More recently, NS-Leibniz bialgebras
were introduced in~\cite{LW26},
where matched pairs, Manin triples,
Yang-Baxter equations and relative
Rota-Baxter operators were studied.
These results provide one of the main
motivations for developing a bialgebra
theory for NS-algebras.

The construction of a bialgebra theory
for NS-algebras requires an appropriate
description of the compatibility between
their three algebra operations and the
corresponding coalgebra operations.
A natural approach is to investigate 
representations, dual representations,
matched pairs and Manin triples of NS-algebras.
In particular, the passage from a
representation to its dual requires
careful treatment, since the three
operations of an NS-algebra are subject
to different compatibility identities.
This motivates the study of the
coadjoint representation and its role
in constructing an NS-algebra structure
on the direct sum of an NS-algebra
and its dual.

The operator approach to Yang-Baxter equations
provides a framework for constructing solutions
via algebraic operators.
Bai~\cite{Bai07} developed a unified algebraic
approach to the classical Yang-Baxter equation,
establishing connections with $\mathcal O$-operators
and pre-Lie algebras.
Motivated by this approach, we investigate
the relationship between symmetric solutions
of the NS-Yang-Baxter equation and
$\mathcal O$-operators on NS-algebras.
We further introduce pre-NS-algebras
and study their relationship with
invertible $\mathcal O$-operators.

A related development concerns dendriform
D-bialgebras and Connes cocycles.
Bai~\cite{Bai1} established correspondences
between double constructions of Frobenius
algebras and antisymmetric infinitesimal
bialgebras, and between double constructions
of Connes cocycles and dendriform
D-bialgebras.
The relationship between Nijenhuis operators
and associative D-bialgebras was investigated
in~\cite{ML24}.
Dendriform-Nijenhuis bialgebras and
the associated DN-associative
Yang-Baxter equation were
investigated in~\cite{PZGL21}.

Further developments, including
representations and quasi-triangular
structures of dendriform-Nijenhuis
bialgebras, were obtained in~\cite{LM26}.
More recently, Wang~\cite{Wang25}
developed the theory of quasi-triangular
and factorizable dendriform D-bialgebras,
establishing connections with relative
Rota-Baxter operators and quadratic
Rota-Baxter dendriform algebras.
Hou~\cite{Hou27} investigated the
affinization of dendriform D-bialgebras
and its connections with Lie bialgebras
and solutions of the classical
Yang-Baxter equation.
Building on the correspondence between
double constructions of Connes cocycles
and dendriform D-bialgebras,
we establish a construction of
dendriform D-bialgebras from NS-bialgebras.

The organization of the paper and its main results
are as follows.
In Section~\mref{sec:nsa}, we develop the
representation theory of NS-algebras.
We characterize representations in terms of
semidirect products, construct their duals,
and investigate the relationship between
the adjoint and coadjoint representations
of quadratic NS-algebras.

In Section~\mref{sec:nsb}, we introduce
matched pairs, Manin triples and NS-bialgebras,
and establish their equivalence
(Theorem~\mref{thm:equiv-ns}).
In particular, an NS-bialgebra structure
on $A$ is characterized by a matched pair
of NS-algebras on $A$ and $A^*$, or equivalently,
by a Manin triple structure on $A\bowtie A^*$.

In Section~\mref{sec:quasi-triangular-ns}, we introduce
the NS-YBE and
investigate its relationship with triangular
NS-bialgebras and $\mathcal O$-operators.
We prove that every symmetric solution of
the NS-YBE induces a triangular NS-bialgebra
(Proposition~\mref{pro:6.2}). We further establish a correspondence
between $\mathcal O$-operators and symmetric
solutions of the NS-YBE in suitable semidirect
product NS-algebras
(Theorem~\mref{thm:semi-dybe}).
We also introduce pre-NS-algebras, characterize
them in terms of invertible $\mathcal O$-operators,
and construct nondegenerate symmetric solutions
of the NS-YBE from pre-NS-algebras
(Theorem~\mref{thm:pre-ns}).

In Section~\mref{sec:Dbia}, we investigate
applications of NS-bialgebras to dendriform
D-bialgebras. We show that a quadratic
NS-algebra induces a nondegenerate Connes
cocycle on an associated associative algebra.
Using the Manin triple associated with an
NS-bialgebra, we obtain a double construction
of a Connes cocycle. Consequently, we prove
that every NS-bialgebra on $A$ induces a
dendriform D-bialgebra structure on
$A\oplus A$ (Theorem~\mref{thm:ns-to-dendriform}).

\noindent

Throughout this paper, $\bfk$ denotes a field
of characteristic zero. Unless otherwise
stated, all vector spaces and algebras are
finite-dimensional, and all linear maps
and tensor products are over $\bfk$.

\section{NS-algebras and their representations}
\mlabel{sec:nsa}

In this section, we recall NS-algebras and develop
their representation theory. We characterize
representations in terms of semidirect products,
construct their dual representation, and investigate the
relationship between the adjoint and coadjoint
representations of quadratic NS-algebras.
\begin{defi}\cite{Ler04}\mlabel{re:nsa}
	An {\bf NS-algebra} is a quadruple
	$(A,\prec,\succ,\bullet)$, where $A$ is a vector space
	equipped with three bilinear operations
	$\prec,\succ,\bullet:A\otimes A\to A$ satisfying
	the following identities:
	\begin{align}
		(x\prec y)\prec z
		&=x\prec(y\star z),\mlabel{eq:nsa01}\\
		(x\succ y)\prec z
		&=x\succ(y\prec z),\mlabel{eq:nsa02}\\
		(x\star y)\succ z
		&=x\succ(y\succ z),\mlabel{eq:nsa03}\\
		(x\star y)\bullet z+(x\bullet y)\prec z
		&=x\succ(y\bullet z)+x\bullet(y\star z),\quad x,y,z\in A,
		\mlabel{eq:nsa04}
	\end{align}
	where 
	\begin{align}
		x\star y:=x\prec y+x\succ y+x\bullet y. \mlabel{eq:ast}
	\end{align}
\end{defi}

\begin{rmk}\mlabel{rmk:da}
According to \cite{JL01}, a {\bf dendriform algebra} is a triple $(A,\prec,\succ)$, where $A$ is a vector space and $\prec,\succ$ are binary operations on $A$ satisfying
the following identities:
\begin{align*}
	(x\prec y)\prec z
	&=x\prec(y\ast z),\\
	(x\succ y)\prec z
	&=x\succ(y\prec z),\\
	(x\ast y)\succ z
	&=x\succ(y\succ z),\quad x,y,z\in A
\end{align*}
where 
\begin{align}
	x\ast y:=x\prec y+x\succ y. \mlabel{eq:astoda}
\end{align}
Note that an NS-algebra \((A,\prec,\succ,\bullet)\) with \(\bullet=0\) is precisely a dendriform algebra. Consequently, NS-algebras can be viewed as a generalization of dendriform algebras.
	\end{rmk}

\begin{ex}\mlabel{ex:NSa}
	Let $A=\RR^2$ with standard basis
	$\{e_1,e_2\}$. Define three binary
	operations $\prec,\succ,\bullet$ on $A$ by
	\[
	e_1\prec e_1=e_2,\qquad
	e_1\succ e_1=e_2,\qquad
	e_1\bullet e_1=-e_2,
	\]
	with all other products of basis elements
	equal to zero.
	Then $(A,\prec,\succ,\bullet)$ is an
	NS-algebra, since all iterated products
	of three elements vanish.
\end{ex}

\begin{defi}\mlabel{defi:NSAI}
	Let $(A_1,\prec_1,\succ_1,\bullet_1)$ and $(A_2,\prec_2,\succ_2,\bullet_2)$ be  NS-algebras. They are called
	{\bf isomorphic} if there exists a linear isomorphism
	$f:A_1\to A_2$ such that, for all $x,y\in A_1$, 
	\begin{align*}
		f(x\circ_1 y)=f(x)\circ_2 f(y), \quad \circ\in\{\prec,\succ,\bullet\}.
	\end{align*}
\end{defi}

\begin{defi}\mlabel{re:rotnsa}
	Let $(A,\prec,\succ,\bullet)$ be an NS-algebra
	and let $V$ be a vector space.
	A $7$-tuple
	$(V,l_{\prec},
	r_{\prec}$,\\$l_{\succ},r_{\succ},
	l_{\bullet},r_{\bullet})$
	is called a {\bf representation} of
	$(A,\prec,\succ,\bullet)$ if the six linear maps
	$l_{\prec},r_{\prec},l_{\succ},r_{\succ},
	l_{\bullet},r_{\bullet}:A\to\End(V)$
	satisfy the following identities, for all $x,y\in A$:
	\begin{align}
		l_{\prec}(x\prec y)
		&=l_{\prec}(x)l_{\star}(y),
		\mlabel{eq:ronsa01}\\
		r_{\prec}(x)l_{\prec}(y)
		&=l_{\prec}(y)r_{\star}(x),
		\mlabel{eq:ronsa02}\\
		r_{\prec}(x)r_{\prec}(y)
		&=r_{\prec}(y\star x),
		\mlabel{eq:ronsa03}\\
		l_{\prec}(x\succ y)
		&=l_{\succ}(x)l_{\prec}(y),
		\mlabel{eq:ronsa04}\\
		r_{\prec}(x)l_{\succ}(y)
		&=l_{\succ}(y)r_{\prec}(x),
		\mlabel{eq:ronsa05}\\
		r_{\succ}(x\prec y)
		&=r_{\prec}(y)r_{\succ}(x),
		\mlabel{eq:ronsa06}\\
		r_{\succ}(x\succ y)
		&=r_{\succ}(y)r_{\star}(x),
		\mlabel{eq:ronsa07}\\
		l_{\succ}(x)r_{\succ}(y)
		&=r_{\succ}(y)l_{\star}(x),
		\mlabel{eq:ronsa08}\\
		l_{\succ}(x)l_{\succ}(y)
		&=l_{\succ}(x\star y),
		\mlabel{eq:ronsa09}\\
		l_{\bullet}(x\star y)+l_{\prec}(x\bullet y)
		&=l_{\succ}(x)l_{\bullet}(y)
		+l_{\bullet}(x)l_{\star}(y),
		\mlabel{eq:ronsa10}\\
		l_{\succ}(x)r_{\bullet}(y)
		+l_{\bullet}(x)r_{\star}(y)
		&=r_{\prec}(y)l_{\bullet}(x)
		+r_{\bullet}(y)l_{\star}(x),
		\mlabel{eq:ronsa11}\\
		r_{\bullet}(x\star y)+r_{\succ}(x\bullet y)
		&=r_{\prec}(y)r_{\bullet}(x)
		+r_{\bullet}(y)r_{\star}(x).
		\mlabel{eq:ronsa12}
	\end{align}
	Here $\star$ is defined by Eq.~\meqref{eq:ast}, and
	\begin{equation}
		l_{\star}:=l_{\prec}+l_{\succ}+l_{\bullet},
		\qquad
		r_{\star}:=r_{\prec}+r_{\succ}+r_{\bullet}.
		\mlabel{eq:lrast}
	\end{equation}
\end{defi}

\begin{pro}\mlabel{pro:spnsa}
	Let $(A,\prec,\succ,\bullet)$ be an NS-algebra
	and let $V$ be a vector space.
	Let
	$l_{\prec},r_{\prec},l_{\succ},r_{\succ},
	l_{\bullet},\\r_{\bullet}:A\to\End(V)$
	be six linear maps. Define three bilinear operations
	$\prec_{\ltimes},\succ_{\ltimes},\bullet_{\ltimes}:
	(A\oplus V)\otimes(A\oplus V)\to A\oplus V$
	by
	\begin{align}
		(x,u)\circ_{\ltimes}(y,v)
		&:=(x\circ y,l_{\circ}(x)v+r_{\circ}(y)u), \quad  \circ\in \{\prec,\succ,\bullet\},\mlabel{eq:spnsa}
	\end{align}
	for all $(x,u),(y,v)\in A\oplus V$.
	Then
	$(V,l_{\prec},r_{\prec},l_{\succ},r_{\succ},
	l_{\bullet},r_{\bullet})$
	is a representation of $(A,\prec,\succ,\bullet)$
	if and only if
	$(A\oplus V,\prec_{\ltimes},\succ_{\ltimes},
	\bullet_{\ltimes})$
	is an NS-algebra. In this case, the latter
	is called the {\bf semidirect product NS-algebra}
	of $A$ and $V$, denoted by $A\ltimes^{l_{\prec},l_{\succ},l_{\bullet}}_{r_{\prec},r_{\succ},r_{\bullet}} V$, or simply $A\ltimes V$.
\end{pro}

\begin{proof}
	For all $(x,u),(y,v),(z,w)\in A\oplus V$, we have 
	\begin{align*}
	\quad((x,u)\prec_{\ltimes} (y,v))\prec_{\ltimes} (z,w)&=(x\prec y,l_{\prec}(x)v+r_{\prec}(y)u)\prec_{\ltimes} (z,w)\\
	&=((x\prec y)\prec z,l_{\prec}(x\prec y)w+r_{\prec}(z)(l_{\prec}(x)v+r_{\prec}(y)u)),
	\end{align*}and\begin{align*}
&(x,u)\prec_{\ltimes}((y,v)\star_{\ltimes} (z,w))
	=(x,u)\prec_{\ltimes}(y\prec z,l_{\prec}(y)w+r_{\prec}(z)v)\\
	&\quad+(x,u)\prec_{\ltimes}(y\succ z,l_{\succ}(y)w+r_{\succ}(z)v)
	+(x,u)\prec_{\ltimes}(y\bullet z,l_{\bullet}(y)w+r_{\bullet}(z)v)\\
	&=(x\prec (y\prec z),l_{\prec}(x)(l_{\prec}(y)w+r_{\prec}(z)v)+r_{\prec}(y\prec z)u)\\
	&\quad
	+(x\prec (y\succ z),l_{\prec}(x)(l_{\succ}(y)w+r_{\succ}(z)v)+r_{\prec}(y\succ z)u)\\
	&\quad
	+(x\prec (y\bullet z),l_{\prec}(x)(l_{\bullet}(y)w+r_{\bullet}(z)v)+r_{\prec}(y\bullet z)u),
\end{align*}	 
	which shows that $(A\oplus V,\prec_{\ltimes},\succ_{\ltimes},
	\bullet_{\ltimes})$ satisfies Eq.~(\mref{eq:nsa01}) if and only if Eqs.~(\mref{eq:ronsa01})--(\mref{eq:ronsa03}) hold. Similarly,  $(A\oplus V,\prec_{\ltimes},\succ_{\ltimes},
	\bullet_{\ltimes})$ satisfies Eqs.~(\mref{eq:nsa02}),~(\mref{eq:nsa03}) and~(\mref{eq:nsa04}) if and only if Eqs.~(\mref{eq:ronsa04})--(\mref{eq:ronsa06}),~(\mref{eq:ronsa07})--(\mref{eq:ronsa09}) and~(\mref{eq:ronsa10})--(\mref{eq:ronsa12}) hold, respectively.
	Thus, $(V,l_{\prec},r_{\prec},l_{\succ},r_{\succ},
	l_{\bullet},r_{\bullet})$
	is a representation of $(A,\prec,\succ,\bullet)$
	if and only if
	$(A\oplus V,\prec_{\ltimes},\succ_{\ltimes},
	\bullet_{\ltimes})$
	is an NS-algebra.
\end{proof}

\begin{pro}
	Let $(A,\prec,\succ,\bullet)$ be an NS-algebra.
	Then
	$(A,L_{\prec},R_{\prec},L_{\succ},R_{\succ},
	L_{\bullet},R_{\bullet})$
	is a representation of $A$, called the
	{\bf adjoint representation} of $A$.
\end{pro}

\begin{lem}\mlabel{lem:nsaaa}
	Let $A$ be a vector space, and let $\prec , \succ ,\bullet: A \ot A\to A$ be three binary operations. Then $(A, \prec , \succ ,\bullet )$ is an NS-algebra if and only if $(A,\star)$ is an associative algebra and $(A,L_{\succ},R_{\prec})$ is a representation of $(A,\star)$, where the operation $\star$ is given by Eq.~\meqref{eq:ast}.
	In this case, $(A,\star)$ is called the
	{\bf associated associative algebra} of
	$(A,\prec,\succ,\bullet)$, and
	$(A,\prec,\succ,\bullet)$ is called a
	{\bf compatible NS-algebra} of $(A,\star)$.
\end{lem}
\begin{proof}
	($\Rightarrow$) Suppose that $(A,\prec,\succ,\bullet)$
	is an NS-algebra. By~\cite[Definition~2.1]{Ler04}, $(A,\star)$ is an associative algebra. For all $x,y,z\in A$,
	we have
	\begin{align*}
		 L_{\succ}(x\star y)z
		&=(x\star y)\succ z
		\stackrel{(\mref{eq:nsa03})}=x\succ(y\succ z)
		=L_{\succ}(x)L_{\succ}(y)z,\\
		R_{\prec}(x\star y)z
		&=z\prec(x\star y)
		\stackrel{(\mref{eq:nsa01})}=(z\prec x)\prec y
		=R_{\prec}(y)R_{\prec}(x)z,\\
		L_{\succ}(x)R_{\prec}(y)z
		&=x\succ(z\prec y)
		\stackrel{(\mref{eq:nsa02})}=(x\succ z)\prec y
		=R_{\prec}(y)L_{\succ}(x)z.
	\end{align*}
	Thus, $(A,L_{\succ},R_{\prec})$ is a
	representation of $(A,\star)$.

$(\Leftarrow)$ Suppose that $(A,\star)$
is an associative algebra and
$(A,L_\succ,R_\prec)$ is a
representation of $(A,\star)$.
The representation identities imply
Eqs.~\meqref{eq:nsa01}--\meqref{eq:nsa03}.
Moreover, using these three identities,
we obtain
\begin{align}
	(x\star y)\bullet z
	+(x\bullet y)\prec z
-x\succ(y\bullet z)
	-x\bullet(y\star z)
	=(x\star y)\star z-x\star(y\star z)
	=0.\mlabel{eq:po4}
\end{align}
Thus Eq.~\meqref{eq:nsa04} also holds,
and $(A,\prec,\succ,\bullet)$
is an NS-algebra.
\end{proof}

\begin{pro}\mlabel{pro:roaa}
	Let $(A,\prec,\succ,\bullet)$ be an NS-algebra, and let $(V,l_{\prec},r_{\prec},l_{\succ},r_{\succ},
	l_{\bullet},r_{\bullet})$ be a representation of $A$. Then $(A,\star)$ is an associative algebra, and $(V,l_{\star},r_{\star})$ is a representation of $(A,\star)$, where 
	the operation $\star $ and maps
	$l_{\star },r_{\star }$ are defined by Eqs.~\meqref{eq:ast} and \meqref{eq:lrast}, respectively.
\end{pro}
\begin{proof}
	By Proposition~\mref{pro:spnsa},
	$A\ltimes V$ is an NS-algebra.
	Hence, by Lemma~\mref{lem:nsaaa},
	its associated algebra
	$(A\oplus V,\star_{\ltimes})$
	is associative, where
	\[
	(x,u)\star_{\ltimes}(y,v)
	=(x\star y,l_\star(x)v+r_\star(y)u).
	\]
	The associativity of $\star_{\ltimes}$
	implies that $(A,\star)$ is associative
	and that $(V,l_\star,r_\star)$
	is a representation of $(A,\star)$.
\end{proof}

\begin{cor}
	Let $(A,\prec_A,\succ_A,\bullet_A)$ be an NS-algebra, and let $(D,\prec_D,\succ_D)$ be a dendriform algebra. Then $(A\ot D,\prec_{\ot},\succ_{\ot},\bullet_{\ot})$ is an NS-algebra,  where $\prec_{\ot},\succ_{\ot},\bullet_{\ot}$ are defined by
	\begin{align*}
		(a\ot x)\prec_{\ot}(b\ot y):&=(a\prec_{A} b)\ot(x\prec_{D} y),\\
		(a\ot x)\succ_{\ot}(b\ot y):&=(a\succ_{A} b)\ot(x\succ_{D} y),\\
		(a\ot x)\bullet_{\ot}(b\ot y):&=(a\prec_{A} b)\ot(x\succ_{D} y)+(a\succ_{A} b)\ot(x\prec_{D} y)+(a\bullet_{A} b)\ot(x\ast_{D} y),
	\end{align*} 
for all $a\ot x,b\ot y\in A\ot D$, $\star_{\ot}:=\prec_{\ot}+\succ_{\ot}+\bullet_{\ot}$, and $\star_{A},\ast_{D}$ are given by Eq.~\meqref{eq:ast} and Eq.~\meqref{eq:astoda}, respectively. 
\end{cor}
\begin{proof}
Since $(D,\prec_D,\succ_D)$ is a dendriform algebra, by \cite{JL01}, $(D,\ast_D)$ is an associative algebra. Since $(A,\prec_A,\succ_A,\bullet_A)$ is an NS-algebra, by Lemma~\mref{lem:nsaaa},  $(A,\star_{A})$ is an associative algebra. Thus, $(A\ot D,\star_{\ot})$ is an associative algebra. For all $a\ot x,b\ot y,c\ot z\in A\ot D$, we have
 \begin{align*}
((a\ot x)\prec_{\ot}(b\ot y))\prec_{\ot}(c\ot z)
&=((a\prec_{A} b)\ot(x\prec_{D} y))\prec_{\ot}(c\ot z)\\
&=((a\prec_{A} b)\prec_{A}c)\ot((x\prec_{D} y)\prec_{D}z)\\
&=(a\prec_{A}(b\star_{A} c))\ot(x\prec_{D}(y\ast_{D} z))\\
&=(a\ot x)\prec_{\ot}((b\star_{A} c)\ot(y\ast_{D} z))\\
&=(a\ot x)\prec_{\ot}((b\ot y)\star_{\ot}(c\ot z)).
\end{align*}
which shows that $(A\ot D,\prec_{\ot},\succ_{\ot},\bullet_{\ot})$ satisfies Eq.~\meqref{eq:nsa01}. Similarly, $(A\ot D,\prec_{\ot},\succ_{\ot},\bullet_{\ot})$ satisfies Eqs.~\meqref{eq:nsa02}--\meqref{eq:nsa03}, and then, similar to Eq.~\meqref{eq:po4}, $(A\ot D,\prec_{\ot},\succ_{\ot},\bullet_{\ot})$ satisfies Eq.~\meqref{eq:nsa04}. Therefore, $(A\ot D,\prec_{\ot},\succ_{\ot},\bullet_{\ot})$ is an NS-algebra. 
\end{proof}

In the sequel, we fix the following notations. For a vector space \(V\), we denote its dual space by \(V^*\). The natural pairing between \(V^*\) and \(V\) is given by
\[
\langle\,,\,\rangle : V^*\times V\to \bfk,
\qquad
\langle u^*, v \rangle := u^*(v),
\quad u^*\in V^*,\ v\in V.
\]
Let $W$ be a vector space. For a linear map \(\psi:V\to W\), its linear dual is defined by
\begin{equation*}
	\psi^*: W^*\to V^*,\qquad
	\langle \psi^*(w^*), v \rangle := \langle w^*, \psi(v) \rangle,
	\quad v\in V,\ w^*\in W^*.
\end{equation*}

Let \(A\) be a vector space and \(\rho:A\to\mathrm{End}(V)\) be a linear map. We define the linear map \(\rho^{*}:A\to\mathrm{End}(V^{*})\) by
\begin{equation*}
	\langle \rho^{*}(x)u^{*},v\rangle
	=\langle u^{*},\rho(x)v\rangle,
	\quad x\in A,\ u^{*}\in V^{*},\ v\in V.
\end{equation*}

\begin{pro}\mlabel{pro:dul-rep}
	Let $(A,\prec,\succ,\bullet)$ be an NS-algebra,
	and let
	$(V,l_{\prec},r_{\prec},l_{\succ},r_{\succ},
	l_{\bullet},r_{\bullet})$
	be a representation of $(A,\prec,\succ,\bullet)$.
	Then
	$
	(V^*,-r_{\succ}^*,l_{\star}^*,
	r_{\star}^*,-l_{\prec}^*,
	-r_{\bullet}^*,-l_{\bullet}^*)
	$
	is also a representation of
	$(A,\prec,\succ,\bullet)$,
	where 
	the operation $\star $ and maps
	$l_{\star },r_{\star }$ are defined by Eqs.~\meqref{eq:ast} and \meqref{eq:lrast}, respectively.
\end{pro}
\begin{proof}
Since $(A,\prec,\succ,\bullet)$ is an NS-algebra and $(V,l_{\prec},r_{\prec},
	l_{\succ},r_{\succ},l_{\bullet},r_{\bullet})$
	is a representation of $A$, by Proposition \mref{pro:roaa}, $(V,l_{\star},r_{\star})$ is a representation of $(A,\star)$.
	Thus, for all $x,y\in A$, we have
	\begin{align}
		l_{\star}(x\star y)
		=l_{\star}(x)l_{\star}(y),
		\quad
		r_{\star}(x\star y)
		=r_{\star}(y)r_{\star}(x),
		\quad
		l_{\star}(x)r_{\star}(y)
		=r_{\star}(y)l_{\star}(x).
		\label{eq:dualassoc}
	\end{align}
	For all $x,y\in A$, $u\in V$ and $\xi\in V^*$, we have
	\begin{align*}
		\langle -r_{\succ}^*(x\prec y)\xi,u\rangle=\langle \xi,-r_{\succ}(x\prec y)u\rangle
		\stackrel{(\mref{eq:ronsa06})}=\langle \xi,-r_{\prec}(y)r_{\succ}(x)u\rangle
		=\langle -r_{\succ}^*(x)r^*_{\prec}(y)\xi,u\rangle,
	\end{align*}
	which shows that $(V^*,-r_{\succ}^*,l_{\star}^*,
	r_{\star}^*,-l_{\prec}^*,
	-r_{\bullet}^*,-l_{\bullet}^*)$ satisfies Eq.~(\mref{eq:ronsa01}).
	Similarly, by Eqs.~(\mref{eq:ronsa01})--(\mref{eq:ronsa12})  and~(\mref{eq:dualassoc}), $(V^*,-r_{\succ}^*,l_{\star}^*,
	r_{\star}^*,-l_{\prec}^*,
	-r_{\bullet}^*,-l_{\bullet}^*)$ satisfies  Eqs.~(\mref{eq:ronsa02})--(\mref{eq:ronsa12}). Therefore, $
	(V^*,-r_{\succ}^*,l_{\star}^*,
	r_{\star}^*,-l_{\prec}^*,
	-r_{\bullet}^*,\\-l_{\bullet}^*)$ 
	is a representation of
	$(A,\prec,\succ,\bullet)$.
\end{proof}

\begin{cor}
	Let $(A,\prec,\succ,\bullet)$ be an NS-algebra.
	Then
$
	(A^*,-R_{\succ}^*,L_{\star}^*,
	R_{\star}^*,-L_{\prec}^*,
	-R_{\bullet}^*,-L_{\bullet}^*)
$
	is a representation of $(A,\prec,\succ,\bullet)$,
	called the {\bf coadjoint representation} of $(A,\prec,\succ,\bullet)$,
	where $\star$ is defined by Eq.~\meqref{eq:ast}.
\end{cor}

\begin{proof}
	This follows immediately from
	Proposition~\mref{pro:dul-rep} by taking
	the adjoint representation of $(A,\prec,\succ,\bullet)$.
\end{proof}

\begin{defi}
	Let $(A,\prec,\succ,\bullet)$ be an NS-algebra,
	and let $\frakB:A\times A\to\bfk$ be a bilinear form.
	\begin{enumerate}
		\item The bilinear form $\frakB$ is called
		{\bf nondegenerate} if
		$\frakB(x,y)=0$ for all $y\in A$ implies $x=0$.
		
		\item The bilinear form $\frakB$ is called
		{\bf antisymmetric} if
		\[
		\frakB(x,y)=-\frakB(y,x),
		\qquad x,y\in A.
		\]
		
		\item The bilinear form $\frakB$ is called
		{\bf invariant} if, for all $x,y,z\in A$,
		\begin{align}
			\frakB(x\prec y,z)
			&=\frakB(x,y\star z),
			\mlabel{eq:ibf2}\\
			\frakB(x\succ y,z)
			&=\frakB(y,z\star x),
			\mlabel{eq:ibf1}\\
			\frakB(x\bullet y,z)
			&=-\frakB(x,y\bullet z),
			\mlabel{eq:ibf3}
		\end{align}
		where $\star$ is defined by
		Eq.~\meqref{eq:ast}.
	\end{enumerate}
\end{defi}

\begin{pro}
	Let $(A,\prec,\succ,\bullet)$ be an NS-algebra,
	and let $\frakB:A\times A\to\bfk$ be an
	antisymmetric invariant bilinear form.
	Then, for all $x,y,z\in A$,
	\begin{equation*}
		\begin{aligned}
			&\frakB(x\star y,z)
			+\frakB(y\star z,x)
			+\frakB(z\star x,y)=\frakB(x\bullet y,z).
		\end{aligned}
	\end{equation*}
	In particular, if $\frakB$ is nondegenerate,
	then the following statements are equivalent:
	\begin{enumerate}
		\item
		$\frakB(x\star y,z)+\frakB(y\star z,x)
		+\frakB(z\star x,y)=0$ 
		for all $x,y,z\in A$.
		\item
		$x\bullet y=0$ for all $x,y\in A$.
	\end{enumerate}
\end{pro}

\begin{proof}
	By the invariance and antisymmetry of
	$\frakB$, we have
	\begin{align*}
		\frakB(x\star y,z)
		&=\frakB(x\prec y,z)
		+\frakB(x\succ y,z)
		+\frakB(x\bullet y,z)\\
		&=\frakB(x,y\star z)
		+\frakB(y,z\star x)
		+\frakB(x\bullet y,z)\\
		&=-\frakB(y\star z,x)
		-\frakB(z\star x,y)
		+\frakB(x\bullet y,z).
	\end{align*}
	This proves the asserted identity.
	The equivalence follows immediately
	from the nondegeneracy of $\frakB$.
\end{proof}

\begin{defi}
	A {\bf quadratic NS-algebra} is a quintuple
	$(A,\prec,\succ,\bullet,\frakB)$,
	where $(A,\prec,\succ,\bullet)$ is an NS-algebra
	and $\frakB:A\times A\to\bfk$ is a
	nondegenerate antisymmetric invariant
	bilinear form.
\end{defi}

\begin{defi}\mlabel{defi:RI}
	Let $(A,\prec,\succ,\bullet)$ be an NS-algebra, and let
	\[
	(V_i,l_{\prec_i},r_{\prec_i},l_{\succ_i},r_{\succ_i},
	l_{\bullet_i},r_{\bullet_i}),\qquad i=1,2,
	\]
	be two representations of $A$. They are called
	{\bf isomorphic} if there exists a linear isomorphism
	$f:V_1\to V_2$ such that, for all $x\in A$,
	\begin{align}
		fl_{\circ_1}(x)=l_{\circ_2}(x)f,\quad
		fr_{\circ_1}(x)=r_{\circ_2}(x)f,\quad \circ\in\{\prec,\succ,\bullet\}.
		\mlabel{eq:RI01}
	\end{align}
\end{defi}

\begin{pro}
	Let $(A,\prec,\succ,\bullet,\frakB)$
	be a quadratic NS-algebra. Then the
	adjoint representation
$
	(A,L_{\prec},R_{\prec},L_{\succ},R_{\succ},
	L_{\bullet},R_{\bullet})
$
	and the coadjoint representation
$
	(A^*,-R_{\succ}^*,L_{\star}^*,
	R_{\star}^*,-L_{\prec}^*,
	-R_{\bullet}^*,-L_{\bullet}^*)
$
	of $(A,\prec,\succ,\bullet)$
	are isomorphic.
\end{pro}

\begin{proof}
	Define a linear map $f:A\to A^*$ by
	$$
	\langle f(x),y\rangle=\frakB(x,y),
	\qquad x,y\in A.
	$$
	Since $\frakB$ is nondegenerate and $A$
	is finite-dimensional, $f$ is a
	linear isomorphism. For $x,y,z\in A$, we have
\begin{align*}
	\langle f(L_{\prec}(x)y),z\rangle&=\frakB(x\prec y,z)\stackrel{\meqref{eq:ibf2}}=\frakB(x,y\star z)\stackrel{\meqref{eq:ibf1}}=\frakB(z\succ x,y )=-\frakB(y,z\succ x)=\langle -R_{\succ}^*(x)f(y),z \rangle,\\
	\langle f(R_{\prec}(x)y),z\rangle&=\frakB(y\prec x,z)\stackrel{\meqref{eq:ibf2}}=\frakB(y,x\star z)=\langle L_{\star}^*(x)f(y),z\rangle,\\
	\langle f(L_{\succ}(x)y),z\rangle&=\frakB(x\succ y,z)\stackrel{\meqref{eq:ibf1}}=\frakB(y,z\star x)=\langle R_{\star}^*(x)f(y),z\rangle,\\
	\langle f(R_{\succ}(x)y),z\rangle&=\frakB(y\succ x,z)\stackrel{\meqref{eq:ibf1}}=\frakB(x,z\star y)\stackrel{\meqref{eq:ibf2}}=\frakB(x\prec z,y)=-\frakB(y,x\prec z)=\langle -L_{\prec}^*(x)f(y),z\rangle,\\
	\langle f(L_{\bullet}(x)y),z\rangle&=\frakB(x\bullet y,z)\stackrel{\meqref{eq:ibf3}}=-\frakB(x,y\bullet z)=\frakB(y\bullet z,x)\stackrel{\meqref{eq:ibf3}}=-\frakB(y,z\bullet x)=\langle -R_{\bullet}^*(x)f(y),z\rangle,\\
	\langle f(R_{\bullet}(x)y),z\rangle&=\frakB(y\bullet x,z)\stackrel{\meqref{eq:ibf3}}=-\frakB(y,x\bullet z)=\langle -L_{\bullet}^*(x)f(y),z\rangle.
\end{align*}
Therefore, by Definition~\mref{defi:RI}, the
adjoint representation
$
(A,L_{\prec},R_{\prec},L_{\succ},R_{\succ},
L_{\bullet},R_{\bullet})
$
and the coadjoint representation
$
(A^*,-R_{\succ}^*,L_{\star}^*,
R_{\star}^*,-L_{\prec}^*,
-R_{\bullet}^*,-L_{\bullet}^*)
$
of $(A,\prec,\succ,\bullet)$
are isomorphic.
\end{proof}

\section{NS-bialgebras, matched pairs and Manin triples of NS-algebras}
\mlabel{sec:nsb}

In this section, we introduce matched pairs and Manin triples of NS-algebras,
and NS-bialgebras. We then characterize matched pairs of
NS-algebras in terms of bicrossed products and investigate
their relationship with Manin triples of NS-algebras. We establish
an equivalence among NS-bialgebras, matched pairs, and
Manin triples of NS-algebras.

\subsection{Matched pairs of NS-algebras}
\begin{defi}\mlabel{def:mpotnsa}
	Let $(A_{1}, \prec_{1}, \succ_{1}, \bullet_{1})$ and $(A_{2}, \prec_{2}, \succ_{2}, \bullet_{2})$ be NS-algebras. For $i,j\in\{1,2\}$ and $i\neq j$, let
	\[
	l_{\prec_i},r_{\prec_i},l_{\succ_i},r_{\succ_i},
	l_{\bullet_i},r_{\bullet_i}:A_i\longrightarrow\End(A_j)
	\]
	be linear maps. The $14$-tuple
	\begin{align*}
		\bigl( (A_1,\prec_1,\succ_1,\bullet_1),
		(A_2,\prec_2,\succ_2,\bullet_2),
		l_{\prec_1},r_{\prec_1},l_{\succ_1},r_{\succ_1},
		l_{\bullet_1},r_{\bullet_1},
		l_{\prec_2},r_{\prec_2},l_{\succ_2},r_{\succ_2},
		l_{\bullet_2},r_{\bullet_2}\bigr)
	\end{align*}
	is called a {\bf matched pair of NS-algebras} if
$
	(A_2,l_{\prec_1},r_{\prec_1},l_{\succ_1},r_{\succ_1},
	l_{\bullet_1},r_{\bullet_1})
$
	and
$
	(A_1,l_{\prec_2},r_{\prec_2},l_{\succ_2},r_{\succ_2},
	l_{\bullet_2},\\r_{\bullet_2})
$
	are representations of $A_1$ and $A_2$, respectively,
	and the following compatibility conditions hold for all
	$x,y\in A_1$ and $u,v\in A_2$:
	\begin{align}
		(l_{\prec_{2}}(u)x)\prec_{1} y+l_{\prec_{2}}(r_{\prec_{1}}(x)u)y&=l_{\prec_{2}}(u)(x\star_{1} y)\mlabel{eq:mponsa01},\\
		x\prec_{1}(r_{\star_{2}}(u)y)+r_{\prec_{2}}(l_{\star_{1}}(y)u)x&=r_{\prec_{2}}(u)(x\prec_{1} y)\mlabel{eq:mponsa02},\\
		(l_{\prec_{1}}(x)u)\prec_{2} v+l_{\prec_{1}}(r_{\prec_{2}}(u)x)v&=l_{\prec_{1}}(x)(u\star_{2} v)\mlabel{eq:mponsa03},\\
		u\prec_{2}(r_{\star_{1}}(x)v)+r_{\prec_{1}}(l_{\star_{2}}(v)x)u&=r_{\prec_{1}}(x)(u\prec_{2} v)\mlabel{eq:mponsa04},\\
		(r_{\prec_{2}}(u)x)\prec_{1} y+l_{\prec_{2}}(l_{\prec_{1}}(x)u)y&=x\prec_{1}(l_{\star_{2}}(u)y)+r_{\prec_{2}}(r_{\star_{1}}(y)u)x\mlabel{eq:mponsa05},\\
		(r_{\prec_{1}}(x)u)\prec_{2} v+l_{\prec_{1}}(l_{\prec_{2}}(u)x)v&=u\prec_{2}(l_{\star_{1}}(x)v)+r_{\prec_{1}}(r_{\star_{2}}(v)x)u\mlabel{eq:mponsa06},\\
		(l_{\succ_{2}}(u)x)\prec_{1} y+l_{\prec_{2}}(r_{\succ_{1}}(x)u)y&=l_{\succ_{2}}(u)(x\prec_{1} y)\mlabel{eq:mponsa07},\\
		x\succ_{1}(r_{\prec_{2}}(u)y)+r_{\succ_{2}}(l_{\prec_{1}}(y)u)x&=r_{\prec_{2}}(u)(x\succ_{1} y)\mlabel{eq:mponsa08},\\
		(l_{\succ_{1}}(x)u)\prec_{2} v+l_{\prec_{1}}(r_{\succ_{2}}(u)x)v&=l_{\succ_{1}}(x)(u\prec_{2} v)\mlabel{eq:mponsa09},\\
		u\succ_{2}(r_{\prec_{1}}(x)v)+r_{\succ_{1}}(l_{\prec_{2}}(v)x)u&=r_{\prec_{1}}(x)(u\succ_{2}v)\mlabel{eq:mponsa10},\\
		(r_{\succ_{2}}(u)x)\prec_{1} y+l_{\prec_{2}}(l_{\succ_{1}}(x)u)y&=x\succ_{1}(l_{\prec_{2}}(u)y)+r_{\succ_{2}}(r_{\prec_{1}}(y)u)x\mlabel{eq:mponsa11},\\
		(r_{\succ_{1}}(x)u)\prec_{2} v+l_{\prec_{1}}(l_{\succ_{2}}(u)x)v&=u\succ_{2}(l_{\prec_{1}}(x)v)+r_{\succ_{1}}(r_{\prec_{2}}(v)x)u\mlabel{eq:mponsa12},\\
		(l_{\star_{2}}(u)x)\succ_{1} y+l_{\succ_{2}}(r_{\star_{1}}(x)u)y&=l_{\succ_{2}}(u)(x\succ_{1} y)\mlabel{eq:mponsa13},\\
		x\succ_{1}(r_{\succ_{2}}(u)y)+r_{\succ_{2}}(l_{\succ_{1}}(y)u)x&=r_{\succ_{2}}(u)(x\star_{1} y)\mlabel{eq:mponsa14},\\
		(l_{\star_{1}}(x)u)\succ_{2} v+l_{\succ_{1}}(r_{\star_{2}}(u)x)v&=l_{\succ_{1}}(x)(u\succ_{2} v)\mlabel{eq:mponsa15},\\
		u\succ_{2}(r_{\succ_{1}}(x)v)+r_{\succ_{1}}(l_{\succ_{2}}(v)x)u&=r_{\succ_{1}}(x)(u\star_{2} v)\mlabel{eq:mponsa16},\\
		(r_{\star_{2}}(u)x)\succ_{1} y+l_{\succ_{2}}(l_{\star_{1}}(x)u)y&=x\succ_{1}(l_{\succ_{2}}(u)y)+r_{\succ_{2}}(r_{\succ_{1}}(y)u)x\mlabel{eq:mponsa17},\\
		(r_{\star_{1}}(x)u)\succ_{2} v+l_{\succ_{1}}(l_{\star_{2}}(u)x)v&=u\succ_{2}(l_{\succ_{1}}(x)v)+r_{\succ_{1}}(r_{\succ_{2}}(v)x)u\mlabel{eq:mponsa18},\\
		l_{\succ_{2}}(u)(x\bullet_{1}y)+l_{\bullet_{2}}(u)(x\star_{1}y)&=(l_{\star_{2}}(u)x)\bullet_{1} y+(l_{\bullet_{2}}(u)x)\prec_{1} y+l_{\bullet_{2}}(r_{\star_{1}}(x)u)y+l_{\prec_{2}}(r_{\bullet_{1}}(x)u)y\mlabel{eq:mponsa19},\\
		r_{\prec_{2}}(u)(x\bullet_{1}y)+r_{\bullet_{2}}(u)(x\star_{1}y)&=x\succ_{1}(r_{\bullet_{2}}(u)y)+x\bullet_{1}(r_{\star_{2}}(u)y)+r_{\succ_{2}}(l_{\bullet_{1}}(y)u)x+r_{\bullet_{2}}(l_{\star_{1}}(y)u)x\mlabel{eq:mponsa20},\\
		l_{\succ_{1}}(x)(u\bullet_{2}v)+l_{\bullet_{1}}(x)(u\star_{2}v)&=(l_{\star_{1}}(x)u)\bullet_{2} v+(l_{\bullet_{1}}(x)u)\prec_{2} v+l_{\bullet_{1}}(r_{\star_{2}}(u)x)v+l_{\prec_{1}}(r_{\bullet_{2}}(u)x)v\mlabel{eq:mponsa21},\\
		r_{\prec_{1}}(x)(u\bullet_{2}v)+r_{\bullet_{1}}(x)(u\star_{2}v)&=u\succ_{2}(r_{\bullet_{1}}(x)v)+u\bullet_{2}(r_{\star_{1}}(x)v)+r_{\succ_{1}}(l_{\bullet_{2}}(v)x)u+r_{\bullet_{1}}(l_{\star_{2}}(v)x)u\mlabel{eq:mponsa22},\\
		(r_{\star_{2}}(u)x)\bullet_{1} y+(r_{\bullet_{2}}(u)x)\prec_{1} y&=r_{\succ_{2}}(r_{\bullet_{1}}(y)u)x+r_{\bullet_{2}}(r_{\star_{1}}(y)u)x+x\succ_{1}(l_{\bullet_{2}}(u)y)\notag\\
		&\quad+x\bullet_{1}(l_{\star_{2}}(u)y)-l_{\bullet_{2}}(l_{\star_{1}}(x)u)y-l_{\prec_{2}}(l_{\bullet_{1}}(x)u)y\mlabel{eq:mponsa23},\\
		(r_{\star_{1}}(x)u)\bullet_{2} v+(r_{\bullet_{1}}(x)u)\prec_{2} v&=r_{\succ_{1}}(r_{\bullet_{2}}(v)x)u+r_{\bullet_{1}}(r_{\star_{2}}(v)x)u+u\succ_{2}(l_{\bullet_{1}}(x)v)\notag\\
		&\quad+u\bullet_{2}(l_{\star_{1}}(x)v)-l_{\bullet_{1}}(l_{\star_{2}}(u)x)v-l_{\prec_{1}}(l_{\bullet_{2}}(u)x)v\mlabel{eq:mponsa24}.
	\end{align}
Here, for $i=1,2$, the operations $\star_i$ and maps
	$l_{\star_i},r_{\star_i}$ are defined by Eqs.~\meqref{eq:ast} and \meqref{eq:lrast}, respectively.
\end{defi}

\begin{pro}\mlabel{pro:matched-ns}
	Let $(A_{1}, \prec_{1}, \succ_{1}, \bullet_{1} )$ and $(A_{2}, \prec_{2}, \succ_{2}, \bullet_{2})$ be NS-algebras, and let 
	\[
	l_{\prec_i},r_{\prec_i},l_{\succ_i},
	r_{\succ_i},l_{\bullet_i},r_{\bullet_i}:
	A_i\to\End(A_j),\qquad i,j\in\{1,2\},\quad i\ne j,
	\]
	be twelve linear maps. Define three bilinear operations
	$\prec_{\bowtie},\succ_{\bowtie},
	\bullet_{\bowtie}$ on $A_1\oplus A_2$ by

\begin{equation}\mlabel{eq:bp}
	\begin{aligned}
		(x,u)\prec_{\bowtie}(y,v)
		&:=\bigl(x\prec_1 y
		+l_{\prec_2}(u)y+r_{\prec_2}(v)x,
		u\prec_2 v
		+l_{\prec_1}(x)v+r_{\prec_1}(y)u\bigr),\\
		(x,u)\succ_{\bowtie}(y,v)
		&:=\bigl(x\succ_1 y
		+l_{\succ_2}(u)y+r_{\succ_2}(v)x,
		u\succ_2 v
		+l_{\succ_1}(x)v+r_{\succ_1}(y)u\bigr),\\
		(x,u)\bullet_{\bowtie}(y,v)
		&:=\bigl(x\bullet_1 y
		+l_{\bullet_2}(u)y+r_{\bullet_2}(v)x,
		u\bullet_2 v
		+l_{\bullet_1}(x)v+r_{\bullet_1}(y)u\bigr),
	\end{aligned}
\end{equation}
	for all $x,y\in A_1,\,u,v\in A_2$. Then the $14$-tuple
$$\bigl( A_1,
		A_2,
		l_{\prec_1},r_{\prec_1},l_{\succ_1},r_{\succ_1},
		l_{\bullet_1},r_{\bullet_1},
		l_{\prec_2},r_{\prec_2},l_{\succ_2},r_{\succ_2},
		l_{\bullet_2},r_{\bullet_2}\bigr)$$ forms a matched pair of
	NS-algebras if and only if
$
	(A_1\oplus A_2,\prec_{\bowtie},
	\succ_{\bowtie},\bullet_{\bowtie})
$
	is an NS-algebra. In this case, it is called
	the {\bf bicrossed product NS-algebra}
	and is denoted by $A_1\bowtie A_2$.
\end{pro}

\begin{proof}
		For all $(x,u),(y,v),(z,w)\in A_1\oplus A_2$, we have 
	\begin{align*}
&\quad((x,u)\prec_{\bowtie} (y,v))\prec_{\bowtie} (z,w)\\
&=(x\prec_1 y
+l_{\prec_2}(u)y+r_{\prec_2}(v)x,
u\prec_2 v
+l_{\prec_1}(x)v+r_{\prec_1}(y)u)\prec_{\bowtie} (z,w)\\
&=((x\prec_1 y
+l_{\prec_2}(u)y+r_{\prec_2}(v)x)\prec_1 z
+l_{\prec_2}(u\prec_2 v
+l_{\prec_1}(x)v+r_{\prec_1}(y)u)z\\&\quad+r_{\prec_2}(w)(x\prec_1 y
+l_{\prec_2}(u)y+r_{\prec_2}(v)x),
(u\prec_2 v
+l_{\prec_1}(x)v+r_{\prec_1}(y)u)\prec_2 w
\\&\quad+l_{\prec_1}(x\prec_1 y
+l_{\prec_2}(u)y+r_{\prec_2}(v)x)w+r_{\prec_1}(z)(u\prec_2 v
+l_{\prec_1}(x)v+r_{\prec_1}(y)u)),\\\\
&\quad(x,u)\prec_{\bowtie}((y,v)\prec_{\bowtie} (z,w))\\
&=(x,u)\prec_{\bowtie}(y\prec_1 z
+l_{\prec_2}(v)z+r_{\prec_2}(w)y,
v\prec_2 w
+l_{\prec_1}(y)w+r_{\prec_1}(z)v )\\
&=(x\prec_1 (y\prec_1 z
+l_{\prec_2}(v)z+r_{\prec_2}(w)y)
+l_{\prec_2}(u)(y\prec_1 z
+l_{\prec_2}(v)z+r_{\prec_2}(w)y)\\&\quad+r_{\prec_2}(v\prec_2 w
+l_{\prec_1}(y)w+r_{\prec_1}(z)v)x,
u\prec_2 (v\prec_2 w
+l_{\prec_1}(y)w+r_{\prec_1}(z)v)
\\&\quad+l_{\prec_1}(x)(v\prec_2 w
+l_{\prec_1}(y)w+r_{\prec_1}(z)v)+r_{\prec_1}(y\prec_1 z
+l_{\prec_2}(v)z+r_{\prec_2}(w)y)u ),\\\\
&\quad(x,u)\prec_{\bowtie}((y,v)\succ_{\bowtie} (z,w))\\
&=(x,u)\prec_{\bowtie}(y\succ_1 z
+l_{\succ_2}(v)z+r_{\succ_2}(w)y,
v\succ_2 w
+l_{\succ_1}(y)w+r_{\succ_1}(z)v )\\
&=(x\prec_1 (y\succ_1 z
+l_{\succ_2}(v)z+r_{\succ_2}(w)y)
+l_{\prec_2}(u)(y\succ_1 z
+l_{\succ_2}(v)z+r_{\succ_2}(w)y)\\&\quad+r_{\prec_2}(v\succ_2 w
+l_{\succ_1}(y)w+r_{\succ_1}(z)v)x,
u\prec_2 (v\succ_2 w
+l_{\succ_1}(y)w+r_{\succ_1}(z)v)
\\&\quad+l_{\prec_1}(x)(v\succ_2 w
+l_{\succ_1}(y)w+r_{\succ_1}(z)v)+r_{\prec_1}(y\succ_1 z
+l_{\succ_2}(v)z+r_{\succ_2}(w)y)u ),
	\end{align*}and\begin{align*}
	&\quad(x,u)\prec_{\bowtie}((y,v)\bullet_{\bowtie} (z,w))\\
	&=(x,u)\prec_{\bowtie}(y\bullet_1 z
	+l_{\bullet_2}(v)z+r_{\bullet_2}(w)y,
	v\bullet_2 w
	+l_{\bullet_1}(y)w+r_{\bullet_1}(z)v )\\
	&=(x\prec_1 (y\bullet_1 z
	+l_{\bullet_2}(v)z+r_{\bullet_2}(w)y)
	+l_{\prec_2}(u)(y\bullet_1 z
	+l_{\bullet_2}(v)z+r_{\bullet_2}(w)y)\\&\quad+r_{\prec_2}(v\bullet_2 w
	+l_{\bullet_1}(y)w+r_{\bullet_1}(z)v)x,
	u\prec_2 (v\bullet_2 w
	+l_{\bullet_1}(y)w+r_{\bullet_1}(z)v)
	\\&\quad+l_{\prec_1}(x)(v\bullet_2 w
	+l_{\bullet_1}(y)w+r_{\bullet_1}(z)v)+r_{\prec_1}(y\bullet_1 z
	+l_{\bullet_2}(v)z+r_{\bullet_2}(w)y)u ),
\end{align*}
which shows that $(A_1\oplus A_2,\prec_{\bowtie},\succ_{\bowtie},\bullet_{\bowtie})$ satisfies Eq.~\meqref{eq:nsa01} if and only if Eqs.~\meqref{eq:mponsa01}--\meqref{eq:mponsa06} hold and $(A_j,l_{\prec_i},r_{\prec_i},l_{\succ_i},
r_{\succ_i},l_{\bullet_i},r_{\bullet_i})$ satisfy Eqs.~\meqref{eq:ronsa01}--\meqref{eq:ronsa03}, where $i,j\in\{1,2\}$ and $i\ne j$. Similarly, $(A_1\oplus A_2,\prec_{\bowtie},\succ_{\bowtie},\bullet_{\bowtie})$ satisfies Eqs.~\meqref{eq:nsa02}--\meqref{eq:nsa04} if and only if Eqs.~\meqref{eq:mponsa07}--\meqref{eq:mponsa24} hold and $(A_j,l_{\prec_i},r_{\prec_i},l_{\succ_i},
r_{\succ_i},l_{\bullet_i},r_{\bullet_i})$ satisfy Eqs.~\meqref{eq:ronsa04}--\meqref{eq:ronsa12}. Thus, $\bigl( A_1,
A_2,
l_{\prec_1},r_{\prec_1},l_{\succ_1},r_{\succ_1},
l_{\bullet_1},r_{\bullet_1},
l_{\prec_2},r_{\prec_2},l_{\succ_2},r_{\succ_2},
l_{\bullet_2},r_{\bullet_2}\bigr)$ is a matched pair of
NS-algebras if and only if
$
(A_1\oplus A_2,\prec_{\bowtie},
\succ_{\bowtie},\bullet_{\bowtie})
$
is an NS-algebra. 	
\end{proof}

\subsection{Manin triples of NS-algebras}

\begin{pro}\mlabel{pro:quadratic-semi-ns}
	Let $(A,\prec,\succ,\bullet)$ be an NS-algebra, and equip $A^*$ with its coadjoint representation
$
	(A,-R^*_{\succ_{A^*}},L_{\star_{A^*}}^*,R^*_{\star_{A^*}},-L^*_{\prec_{A^*}},-R^*_{\bullet_{A^*}},
		-L^*_{\bullet_{A^*}}).
$
	Then the semidirect product $D=A\ltimes A^*$ is a quadratic NS-algebra with respect to the bilinear form $\frakB_d$ defined by
	\begin{equation}
		\frakB_d((x,\xi),(y,\eta))=-\langle x,\eta\rangle+\langle y,\xi\rangle,
		\qquad x,y\in A,\quad \xi,\eta\in A^*.
		\mlabel{eq:Bd1}
	\end{equation}
\end{pro}

\begin{proof}
	Since $(A^*,-R_{\succ_A}^*,L_{\star_A}^*,
	R_{\star_A}^*,-L_{\prec_A}^*,
	-R_{\bullet_A}^*,-L_{\bullet_A}^*)$ is the coadjoint representation of $(A,\prec,\succ,\bullet)$, by Proposition \mref{pro:spnsa}, $A\ltimes A^*$ is an NS-algebra. For all $(x,\xi),(y,\eta),(z,\zeta)\in A\oplus A^*$, we have 
	\begin{align*}
	\frakB_d((x,\xi)\prec_\ltimes (y,\eta),(z,\zeta))
		&=\frakB_d((x\prec y,-R_\succ^*(x)\eta+L_\star^*(y)\xi),(z,\zeta))\\
		&=-\langle x\prec y,\zeta\rangle+\langle z,-R_\succ^*(x)\eta+L_\star^*(y)\xi\rangle\\
		&=-\langle x\prec y,\zeta\rangle-\langle z\succ x,\eta\rangle+\langle y\star z,\xi\rangle.\\
\frakB_d((x,\xi) ,(y,\eta)\star_\ltimes(z,\zeta))
	&=\frakB_d((x,\xi) ,(y\star z, R_\prec^*(y)\zeta+L_\succ^*(z)\eta))\\
	&=-\langle x,R_\prec^*(y)\zeta+L_\succ^*(z)\eta\rangle+\langle y\star z,\xi\rangle\\
	&=-\langle x,R_\prec^*(y)\zeta\rangle-\langle x,L_\succ^*(z)\eta\rangle+\langle y\star z,\xi\rangle\\
	&=-\langle x\prec y,\zeta\rangle-\langle z\succ x,\eta\rangle+\langle y\star z,\xi\rangle,
\end{align*}
	which shows that $(A\ltimes A^*, \frakB_d)$ satisfies Eq.~\meqref{eq:ibf2}.  
	Similarly, the remaining invariance
	identities~Eqs.~\meqref{eq:ibf1}--~\meqref{eq:ibf3} hold. Thus, $(A\ltimes A^*, \frakB_d)$ is a quadratic NS-algebra.
	\end{proof}

\begin{defi}\mlabel{def:manin-ns}
	Let $(A,\prec_A,\succ_A,\bullet_A)$ and $(A^*,\prec_{A^*},\succ_{A^*},\bullet_{A^*})$ be NS-algebras. A {\bf Manin triple of NS-algebras} is a triple
	$
	\bigl((A\oplus A^*,\prec_{\bowtie},\succ_{\bowtie},\bullet_{\bowtie},\frakB_d),A,A^*\bigr)
	$
	such that $(A\oplus A^*,\prec_{\bowtie},\succ_{\bowtie},\bullet_{\bowtie},\frakB_d)$ is a quadratic NS-algebra, $\frakB_d$ is given by Eq.~\meqref{eq:Bd1}, and $A$ and $A^*$ are NS subalgebras.
\end{defi}

\begin{pro}\mlabel{pro:iiandiii}
Let $(A,\prec_A,\succ_A,\bullet_A)$ and $(A^*,\prec_{A^*},\succ_{A^*},\bullet_{A^*})$ be NS-algebras. Then 	$$(A,A^*,-R_{\succ_A}^*,L_{\star_A}^*,
R_{\star_A}^*,-L_{\prec_A}^*,
-R_{\bullet_A}^*,-L_{\bullet_A}^*,-R_{\succ_{A^*}}^*,L_{\star_{A^*}}^*,
	R_{\star_{A^*}}^*,-L_{\prec_{A^*}}^*,
	-R_{\bullet_{A^*}}^*,-L_{\bullet_{A^*}}^*)$$ forms a matched pair of NS-algebras if and only if
	$
	\bigl((A\bowtie A^*,\prec_{\bowtie},\succ_{\bowtie},\bullet_{\bowtie},\frakB_d),A,A^*\bigr)
	$
	is a Manin triple of NS-algebras.
\end{pro}
\begin{proof}
	($\Leftarrow$) Suppose that $
	\bigl((A\bowtie A^*,\prec_{\bowtie},\succ_{\bowtie},\bullet_{\bowtie},\frakB_d),A,A^*\bigr)
	$
	is a Manin triple of NS-algebras. For all $x,y\in A$ and $\xi,\eta\in A^*$, we have
	\begin{align*}
		&\quad\frakB_d((x,0)\prec_{\bowtie}(0,\eta),(y,0))\stackrel{\meqref{eq:ibf2}}=\frakB_d((x,0),(0,\eta)\star_{\bowtie} (y,0))\stackrel{\meqref{eq:ibf1}}=\frakB_d((y,0)\succ_{\bowtie}(x,0),(0,\eta))\\
		&=\frakB_d((y\succ_A x,0) ,(0,\eta))=-\langle R_{\succ_A} (x) y,\eta \rangle=-\langle y,R^*_{\succ_A} (x) \eta \rangle=\frakB_d((0,-R^*_{\succ_A} (x) \eta),(y,0)),\\
	    &\quad\frakB_d((x,0)\prec_{\bowtie}(0,\eta),(0,\xi))\stackrel{\meqref{eq:ibf2}}=\frakB_d((x,0),(0,\eta)\star_{\bowtie}(0,\xi))=\frakB_d((x,0),(0,\eta\star_{A^*}\xi))\\
	    &=-\langle x,L^*_{\star_{A^*}}(\eta)\xi \rangle=-\langle L^*_{\star_{A^*}}(\eta)x,\xi \rangle=\frakB_d(( L^*_{\star_{A^*}}(\eta)x,0),(0, \xi)),
	\end{align*} 
	and 
	\begin{align*}
&\quad\frakB_d((0,\xi)\prec_{\bowtie}(y,0),(x,0))\stackrel{\meqref{eq:ibf2}}=\frakB_d((0,\xi),(y,0)\star_{\bowtie} (x,0))=\frakB_d((0,\xi),(y\star_A x,0) )\\
&=\langle L_\star(y)x,\xi \rangle=\langle x,L^*_{\star_A}(y)\xi \rangle=\frakB_d((0 ,L^*_{\star_A}(y)\xi),(x,0)),\\
&\quad\frakB_d((0,\xi)\prec_{\bowtie}(y,0),( 0,\eta))\stackrel{\meqref{eq:ibf2}}=\frakB_d((0,\xi),(y,0)\star_{\bowtie} (0,\eta))\stackrel{\meqref{eq:ibf1}}=\frakB_d((0,\eta)\succ_{\bowtie} (0,\xi),(y,0))\\
&=\frakB_d((0,\eta\succ_{A^*}\xi) ,(y,0))=\langle y,R_{\succ_{A^*}} (\xi)\eta\rangle=\langle R^*_{\succ_{A^*}}(\xi)y,\eta\rangle=\frakB_d((-R^*_{\succ_{A^*}}(\xi)y,0),( 0,\eta)),
\end{align*}
	which shows 
	$$(x,0)\prec_{\bowtie}(0,\eta)=( L^*_{\star_{A^*}}(\eta)x,-R^*_{\succ_A} (x) \eta),\; (0,\xi)\prec_{\bowtie}(y,0)=(-R_{\succ_{A^*}}^*(\xi)y,L^*_{\star_A}(y)\xi).$$
	Then 
	\begin{align*}
	(x,\xi)\prec_{\bowtie}(y,\eta)&=(x\prec_{A } y,0)+( L^*_{\star_{A^*}}(\eta)x,-R^*_{\succ_A} (x) \eta)+(-R_{\succ_{A^*}}^*(\xi)y,L^*_{\star_{A}}(y)\xi)+(0,\xi\prec_{A^*}\eta)\\
	&=(x\prec_{A } y-R_{\succ_{A^*}}^*(\xi)y+L^*_{\star_{A^*}}(\eta)x,\xi\prec_{A^*}\eta-R^*_{\succ_A} (x) \eta+L^*_{\star_A}(y)\xi).
\end{align*}
	Similarly, we have
	\begin{align*}
		(x,\xi)\succ_{\bowtie}(y,\eta)&=(x\succ_{A} y+R_{\star_{A^*}}^*(\xi)y-L^*_{\prec_{A^*}}(\eta)x,\xi\succ_{A^*}\eta+R^*_{\star_A} (x) \eta-L^*_{\prec_A}(y)\xi),\\
		(x,\xi)\bullet_{\bowtie}(y,\eta)&=(x\bullet_{A} y-R_{\bullet_{A^*}}^*(\xi)y-L^*_{\bullet_{A^*}}(\eta)x,\xi\bullet_{A^*}\eta-R^*_{\bullet_A} (x) \eta-L^*_\bullet(y)\xi).
	\end{align*}
Since $(A^*,-R_{\succ_A}^*,L_{\star_A}^*,
R_{\star_A}^*,-L_{\prec_A}^*,
-R_{\bullet_A}^*,-L_{\bullet_A}^*)$ and $(A,-R_{\succ_{A^*}}^*,L_{\star_{A^*}}^*,
R_{\star_{A^*}}^*,-L_{\prec_{A^*}}^*,
-R_{\bullet_{A^*}}^*,-L_{\bullet_{A^*}}^*)$ are the coadjoint representations of $(A,\prec_A,\succ_A,\bullet_A)$ and $(A^*,\prec_{A^*},\succ_{A^*},\bullet_{A^*})$, respectively, $(A\bowtie A^*,\prec_{\bowtie},\succ_{\bowtie},\bullet_{\bowtie} )$ is an NS-algebra. Thus, by Proposition \mref{pro:matched-ns}, $$(A,A^*,-R_{\succ_A}^*,L_{\star_A}^*,
R_{\star_A}^*,-L_{\prec_A}^*,
-R_{\bullet_A}^*,-L_{\bullet_A}^*,-R_{\succ_{A^*}}^*,L_{\star_{A^*}}^*,
R_{\star_{A^*}}^*,-L_{\prec_{A^*}}^*,
-R_{\bullet_{A^*}}^*,-L_{\bullet_{A^*}}^*)$$ is a matched pair of NS-algebras. 

($\Rightarrow$) Suppose that $$(A,A^*,-R_{\succ_A}^*,L_{\star_A}^*,
	R_{\star_A}^*,-L_{\prec_A}^*,
	-R_{\bullet_A}^*,-L_{\bullet_A}^*,-R_{\succ_{A^*}}^*,L_{\star_{A^*}}^*,
	R_{\star_{A^*}}^*,-L_{\prec_{A^*}}^*,
	-R_{\bullet_{A^*}}^*,-L_{\bullet_{A^*}}^*)$$ is a matched pair of NS-algebras.  Then, by Proposition \mref{pro:matched-ns}, $(A\bowtie A^*,\prec_{\bowtie},\succ_{\bowtie},\bullet_{\bowtie} )$ is an NS-algebra associated with the coadjoint representations, where $\prec_{\bowtie},\succ_{\bowtie},\bullet_{\bowtie}$ are given by Eq.~\meqref{eq:bp}. By Proposition \mref{pro:quadratic-semi-ns}, $(A\bowtie A^*,\prec_{\bowtie},\succ_{\bowtie},\bullet_{\bowtie},\frakB_d)$ is a quadratic NS-algebra, where the bilinear form $\frakB_d$ is given by Eq.~\meqref{eq:Bd1}. For all $(x,0),(y,0)\in A\oplus A^*$, we have
	$$(x,0)\prec_{\bowtie}(y,0)=(x\prec_A y,0),\quad(x,0)\succ_{\bowtie}(y,0)=(x\succ_A y,0),\quad(x,0)\bullet_{\bowtie}(y,0)=(x\bullet_A y,0)\in A.$$
	Thus, $(A,\prec_A,\succ_A,\bullet_A)$ is an NS-subalgebra. Similarly, $(A^*,\prec_{A^*},\succ_{A^*},\bullet_{A^*})$ is also an NS-subalgebra. Therefore, $
	\bigl((A\bowtie A^*,\prec_{\bowtie},\succ_{\bowtie},\bullet_{\bowtie},\frakB_d),A,A^*\bigr)
	$
	is a Manin triple of NS-algebras.
\end{proof}

\subsection{NS-coalgebras and NS-bialgebras}

\begin{defi}~\cite{Das22b}\mlabel{def:ns-coalgebra}
	An {\bf NS-coalgebra}
 is a quadruple $(A,\delta_\prec,\delta_\succ,\delta_\bullet)$ consisting of a vector space $A$ and linear maps $\delta_\prec,\delta_\succ,\delta_\bullet:A\to A\ot A$ satisfying
	\begin{align}
		(\delta_\prec\ot\id)\delta_\prec&=(\id\ot\delta_\star)\delta_\prec,\mlabel{eq:NSca02}\\
		(\delta_\succ\ot\id)\delta_\prec&=(\id\ot\delta_\prec)\delta_\succ,\mlabel{eq:NSca01}\\
		(\delta_\star\ot\id)\delta_\succ&=(\id\ot\delta_\succ)\delta_\succ,\mlabel{eq:NSca03}\\
		(\delta_\star\ot\id)\delta_\bullet+(\delta_\bullet\ot\id)\delta_\prec
		&=(\id\ot\delta_\bullet)\delta_\succ+(\id\ot\delta_\star)\delta_\bullet,\mlabel{eq:NSca04}
	\end{align}
	where
	\begin{equation}
		\delta_\star:=\delta_\prec+\delta_\succ+\delta_\bullet.\mlabel{eq:deltastar}
	\end{equation}
\end{defi}

\begin{rmk}
Let $\prec_{A^*},\succ_{A^*},\bullet_{A^*}:A^*\ot A^*\to A^*$ be bilinear operations. Their  dual linear maps $\delta_{\prec_A},\delta_{\succ_A},\delta_{\bullet_A}:A\to A\ot A$ are uniquely determined by
	\begin{align}
		\langle\xi\circ_{A^*}\eta,x\rangle=\langle\xi\ot\eta,\delta_{\circ_A}(x)\rangle,
		 \quad  \circ\in \{\prec,\succ,\bullet\},\mlabel{eq:lds1}
	\end{align}
	for all $x\in A$ and $\xi,\eta\in A^*$.
\end{rmk}

\begin{pro}\mlabel{pro:dual-ns-coalgebra}
	Let $A$ be a vector space,
	and let $\prec_{A*},\succ_{A*},\bullet_{A*}$
	be three bilinear operations on $A^*$.
	Let $\delta_{\prec_A},\delta_{\succ_A},
	\delta_{\bullet_A}:A\to A\otimes A$
	be their dual linear maps, respectively,
	as defined by Eq.~\meqref{eq:lds1}.
	Then
	$(A,\delta_{\prec_A},\delta_{\succ_A},
	\delta_{\bullet_A})$ is an NS-coalgebra
	if and only if
	$(A^*,\prec_{A*},\succ_{A*},
	\bullet_{A*})$ is an NS-algebra.
\end{pro}
\begin{proof}
	 For all $x\in A$ and $\xi,\eta,\zeta\in A^*$, we have 
	\begin{align*}
	\langle (\xi\prec_{A*} \eta)\prec_{A*} \zeta,x\rangle
    =\langle (\xi\prec_{A*} \eta)\ot \zeta, \delta_{\prec_A}(x)\rangle
	=\langle\xi\ot\eta\ot\zeta,(\delta_{\prec_A}\ot \id)\delta_{\prec_A}(x)\rangle,
	\end{align*} and \begin{eqnarray*} 
	&&\langle \xi\prec_{A*}(\eta\prec_{A*}\zeta+\eta\succ_{A*}\zeta+
	\eta\bullet_{A*} \zeta),x\rangle\\
	&=&\langle\xi\prec_{A*}(\eta\prec_{A*}\zeta),x\rangle+\langle\xi\prec_{A*}(\eta\succ_{A*}\zeta),x\rangle+\langle\xi\prec_{A*}(\eta\bullet_{A*}\zeta),x\rangle\\
	&=&\langle\xi\ot(\eta\prec_{A*}\zeta),\delta_{\prec_A}(x)\rangle+\langle\xi\ot(\eta\succ_{A*}\zeta),\delta_{\prec_A}(x)\rangle+\langle\xi\ot(\eta\bullet_{A*}\zeta),\delta_{\prec_A}(x)\rangle\\
	&=&\langle\xi\ot\eta\ot\zeta,(\id\ot\delta_{\prec_A})\delta_{\prec_A}(x)\rangle+\langle\xi\ot\eta\ot\zeta,(\id\ot\delta_{\succ_A})\delta_{\prec_A}(x)\rangle+\langle\xi\ot\eta\ot\zeta,(\id\ot\delta_{\bullet_A})\delta_{\prec_A}(x)\rangle\\
	&=&\langle\xi\ot\eta\ot\zeta,(\id\ot\delta_{\star_A})\delta_{\prec_A}(x)\rangle,
\end{eqnarray*}
	which shows that $(A^*,\prec_{A*},\succ_{A*},
	\bullet_{A*})$ satisfies Eq.~\meqref{eq:nsa01} if and only if $(A,\delta_{\prec_A},\delta_{\succ_A},
	\delta_{\bullet_A})$ satisfies Eq.~\meqref{eq:NSca02}. Similarly,  we have
	\begin{eqnarray*}
		&\text{Eq}.~\meqref{eq:nsa02}\Leftrightarrow \text{Eq}.~\meqref{eq:NSca01},\quad
		\text{Eq}.~\meqref{eq:nsa03}\Leftrightarrow \text{Eq}.~\meqref{eq:NSca03},\quad
		&\text{Eq}.~\meqref{eq:nsa04}\Leftrightarrow \text{Eq}.~\meqref{eq:NSca04}.
	\end{eqnarray*} Therefore, $(A,\delta_{\prec_A},\delta_{\succ_A},
\delta_{\bullet_A})$ is an NS-coalgebra
if and only if
$(A^*,\prec_{A*},\succ_{A*},
\bullet_{A*})$ is an NS-algebra.
\end{proof}

\begin{ex}\mlabel{ex:NSca}
	Let $A=\RR^2$ with standard basis
	$\{e_1,e_2\}$. Define three linear maps
	$\delta_\prec,\delta_\succ,\delta_\bullet:
	A\to A\ot A$ by
	\[
	\delta_\prec(e_1)=\delta_\succ(e_1)
	=\delta_\bullet(e_1)=0,
	\;
	\delta_\prec(e_2)=\delta_\succ(e_2)
	=-\delta_\bullet(e_2)=e_1\ot e_1.
	\]
	Then $(A,\delta_\prec,\delta_\succ,
	\delta_\bullet)$ is an NS-coalgebra.
	Indeed, all compositions of any two
	of the maps $\delta_\prec,\delta_\succ,
	\delta_\bullet$ vanish.
\end{ex}

\begin{defi}\mlabel{def:ns-bialgebra}
	An {\bf NS-bialgebra} is a septuple $(A,\prec,\succ,\bullet,\delta_\prec,\delta_\succ,\delta_\bullet)$ consisting of an NS-algebra $(A,\prec,\succ,\bullet)$ and an NS-coalgebra $(A,\delta_\prec,\delta_\succ,\delta_\bullet)$ such that, for all $x,y\in A$,
	\begin{align}
		\delta_\prec(x\star y)&=(\id\ot L_\succ(x))\delta_\prec(y)+(R_\star(y)\ot\id)\delta_\prec(x),\mlabel{eq:NSb04}\\
		\delta_\succ(x\star y)&=(R_\prec(y)\ot\id)\delta_\succ(x)+(\id\ot L_\star(x))\delta_\succ(y),\mlabel{eq:NSb01}\\
		\delta_\star(x\prec y)&=(\id\ot L_\prec(x))\delta_\succ(y)+(R_\prec(y)\ot\id)\delta_\star(x),\mlabel{eq:NSb02}\\
		\delta_\star(x\succ y)&=(R_\succ(y)\ot\id)\delta_\prec(x)+(\id\ot L_\succ(x))\delta_\star(y),\mlabel{eq:NSb03}\\
		\big((\id\ot R_\prec(y))-(L_\succ(y)\ot\id)\big)\delta_\star(x)&=\tau\big((L_\prec(x)\ot\id)\delta_\prec(y)-(\id\ot R_\succ(x))\delta_\succ(y)\big),\mlabel{eq:NSb05}\\
		\big((L_\star(x)\ot\id)-(\id\ot R_\prec(x))\big)\delta_\prec(y)&=\tau\big((\id\ot R_\star(y))\delta_\succ(x)-(L_\succ(y)\ot\id)\delta_\succ(x)\big).\mlabel{eq:NSb06}
	\end{align}
	Here $x\star y=x\prec y+x\succ y+x\bullet y$, $\delta_\star$ is defined by Eq.~\meqref{eq:deltastar}, and $\tau:A\ot A\to A\ot A$ is the flip map given by $\tau(a\ot b)=b\ot a$.
\end{defi}

\begin{pro}\mlabel{pro:iandii}
	Let $(A,\prec_A,\succ_A,\bullet_A)$ and
	$(A^*,\prec_{A^*},\succ_{A^*},\bullet_{A^*})$
	be NS-algebras. Let
	$\delta_{\prec_A},\delta_{\succ_A},
	\delta_{\bullet_A}:A\to A\otimes A$
	be the linear duals of
	$\prec_{A^*},\succ_{A^*},\bullet_{A^*}$,
	respectively.
	Then
	$(A,\prec_A,\succ_A,\bullet_A,
	\delta_{\prec_A},\delta_{\succ_A},
	\delta_{\bullet_A})$
	is an NS-bialgebra if and only if
	$$\bigl(A,A^*,
		-R_{\succ_A}^*,L_{\star_A}^*,
		R_{\star_A}^*,-L_{\prec_A}^*,
		-R_{\bullet_A}^*,-L_{\bullet_A}^*,
		-R_{\succ_{A^*}}^*,L_{\star_{A^*}}^*,
		R_{\star_{A^*}}^*,-L_{\prec_{A^*}}^*,
		-R_{\bullet_{A^*}}^*,-L_{\bullet_{A^*}}^*
		\bigr)$$
	is a matched pair of NS-algebras.
\end{pro}
\begin{proof}
	Let $(A,\prec_A,\succ_A,\bullet_A)$ and
	$(A^*,\prec_{A^*},\succ_{A^*},\bullet_{A^*})$
	be NS-algebras. By Proposition \mref{pro:dual-ns-coalgebra}, $(A,\delta_{\prec_A},\delta_{\succ_A},
	\\ \delta_{\bullet_A})$ is an NS-coalgebra
	if and only if
	$(A^*,\prec_{A*},\succ_{A*},
	\bullet_{A*})$ is an NS-algebra. Suppose that, for all $x,y\in A$, 
	\begin{align}
		(R_{\bullet_A}(y)\ot\id)\delta_{\prec_A}(x)+(\id\ot L_{\bullet_A}(x))\delta_{\succ_A}(y)&=(R_{\prec_A}(y)\ot\id)\delta_{\bullet_A}(x)+(\id\ot L_{\succ_A}(x))\delta_{\bullet_A}(y)\notag\\
		&\quad+\delta_{\star_A}(x{\bullet_A} y)-\delta_{\bullet_A}(x{\star_A} y).\mlabel{eq:NSb07}\\
		(\id\ot R_{\bullet_A}(y))\delta_{\succ_A}(x)+\tau(\id\ot R_{\bullet_A}(x))\delta_{\succ_A}(y)&=(L_{\bullet_A}(y)\ot\id)\delta_{\prec_A}(x)+\tau(L_{\bullet_A}(x)\ot\id)\delta_{\prec_A}(y)\notag\\
		&\quad+(\id\ot R_{\prec_A}(y)-L_{\succ_A}(y)\ot\id)\delta_{\bullet_A}(x)\notag\\
		&\quad+\tau(\id\ot R_{\prec_A}(x)-L_{\succ_A}(x)\ot\id)\delta_{\bullet_A}(y).\mlabel{eq:NSb08}
	\end{align} 
	If $(A,\prec_A,\succ_A,\bullet_A,
	\delta_{\prec_A},\delta_{\succ_A},
	\delta_{\bullet_A})$ 
	is an NS-bialgebra, then, by Eqs.~\meqref{eq:NSb04}--\meqref{eq:NSb03} and Eqs.~\meqref{eq:NSb05}--\meqref{eq:NSb06}, we can obtain  Eq.~\meqref{eq:NSb07} and Eq.~\meqref{eq:NSb08}, respectively.
	In Definition~\mref{def:mpotnsa}, take
	$$l_{\prec_1}=-R_{\succ_A}^*,\quad r_{\prec_1}=L_{\star_A}^*,\quad 
	l_{\succ_1}=R_{\star_A}^*,\quad 
	r_{\succ_1}=-L_{\prec_A}^*,\quad l_{\bullet_1}=-R_{\bullet_A}^*,\quad 
	r_{\bullet_1}=-L_{\bullet_A}^*,$$ and
	$$l_{\prec_2}=-R_{\succ_{A^*}}^*,\quad r_{\prec_2}=L_{\star_{A^*}}^*,\quad 
	l_{\succ_2}=R_{\star_{A^*}}^*,\quad 
	r_{\succ_2}=-L_{\prec_{A^*}}^*,\quad l_{\bullet_2}=-R_{\bullet_{A^*}}^*,\quad 
	r_{\bullet_2}=-L_{\bullet_{A^*}}^*.$$
	By Eq.~\meqref{eq:lrast}, we get 
	$$l_{\star_1}=R^*_{\prec_A},
	\quad
	r_{\star_1}=L^*_{\succ_A},
	\quad 
	l_{\star_2}=R_{\prec_{A^*}}^*,
	\quad
	r_{\star_2}=L_{\succ_{A^*}}^*.$$
A routine computation shows that
	\begin{align*}
		\text{\meqref{eq:NSb04}}
		&\Longleftrightarrow
		\text{\meqref{eq:mponsa04},
			\meqref{eq:mponsa09},
			\meqref{eq:mponsa14}},\\
		\text{\meqref{eq:NSb01}}
		&\Longleftrightarrow
		\text{\meqref{eq:mponsa01},
			\meqref{eq:mponsa10},
			\meqref{eq:mponsa15}},\\
		\text{\meqref{eq:NSb02}}
		&\Longleftrightarrow
		\text{\meqref{eq:mponsa02},
			\meqref{eq:mponsa07},
			\meqref{eq:mponsa16}},\\
		\text{\meqref{eq:NSb03}}
		&\Longleftrightarrow
		\text{\meqref{eq:mponsa03},
			\meqref{eq:mponsa08},
			\meqref{eq:mponsa13}},\\
		\text{\meqref{eq:NSb05}}
		&\Longleftrightarrow
		\text{\meqref{eq:mponsa05},
			\meqref{eq:mponsa12},
			\meqref{eq:mponsa17}},\\
		\text{\meqref{eq:NSb06}}
		&\Longleftrightarrow
		\text{\meqref{eq:mponsa06},
			\meqref{eq:mponsa11},
			\meqref{eq:mponsa18}},\\
		\text{\meqref{eq:NSb07}}
		&\Longleftrightarrow
		\text{\meqref{eq:mponsa19}--
			\meqref{eq:mponsa22}},\\
		\text{\meqref{eq:NSb08}}
		&\Longleftrightarrow
		\text{\meqref{eq:mponsa23}--
			\meqref{eq:mponsa24}}.
	\end{align*}
 Eq.~\meqref{eq:NSb07} and Eq.~\meqref{eq:NSb08} follow from Eqs.~\meqref{eq:NSb04}--\meqref{eq:NSb03} and Eqs.~\meqref{eq:NSb05}--\meqref{eq:NSb06}, respectively.
Therefore, 	$(A,\prec_A,\succ_A,\bullet_A,
\delta_{\prec_A},\delta_{\succ_A},
\delta_{\bullet_A})$
is an NS-bialgebra if and only if
$$\bigl(A,A^*,
-R_{\succ_A}^*,L_{\star_A}^*,
R_{\star_A}^*,-L_{\prec_A}^*,
-R_{\bullet_A}^*,-L_{\bullet_A}^*,
-R_{\succ_{A^*}}^*,L_{\star_{A^*}}^*,
R_{\star_{A^*}}^*,-L_{\prec_{A^*}}^*,
-R_{\bullet_{A^*}}^*,-L_{\bullet_{A^*}}^*
\bigr)$$
is a matched pair of NS-algebras.
\end{proof}

\begin{thm}\mlabel{thm:equiv-ns}
	Let $(A,\prec,\succ,\bullet)$ be a
	finite-dimensional NS-algebra, and let
	$(A^*,\prec_*,\succ_*,\bullet_*)$
	be an NS-algebra. Let
	$\delta_{\prec},\delta_{\succ},
	\delta_{\bullet}:A\to A\otimes A$
	be the linear duals of
	$\prec_*,\succ_*,\bullet_*$, respectively.
	
	Then the following statements are equivalent:
	\begin{enumerate}
		\item\mlabel{it:eq1}
		$(A,\prec,\succ,\bullet,
		\delta_{\prec},\delta_{\succ},
		\delta_{\bullet})$ is an NS-bialgebra.
		
		\item\mlabel{it:eq2}
		The coadjoint representations of
		$A$ and $A^*$ form a matched pair
		of NS-algebras.
		
		\item\mlabel{it:eq3}
		The triple
		$
		\bigl((A\bowtie A^*,\frakB_d),A,A^*\bigr)
	$
		is a Manin triple of NS-algebras,
		where $A\bowtie A^*$ is the bicrossed
		product associated with the coadjoint
		representations and $\frakB_d$ is
		defined by Eq.~\meqref{eq:Bd1}.
	\end{enumerate}
\end{thm}

\begin{proof}
	The result follows directly from Propositions~\mref{pro:iiandiii} and~\mref{pro:iandii}. 
\end{proof}

\begin{ex}
	Let $(A,\prec,\succ,\bullet)$ be the
NS-algebra defined in Example~\mref{ex:NSa}.
Then $(A,\prec,\succ,\bullet,
\delta_\prec,\delta_\succ,\\ \delta_\bullet)$
is an NS-bialgebra if and only if
the linear maps $\delta_\prec,\delta_\succ,
\delta_\bullet:A\to A\ot A$
are given by one of the following conditions:
\begin{enumerate}
	\item For all $a,b,c\in\RR$,
	$$\delta_\prec(e_2)=\delta_\succ(e_2)
	=\delta_\bullet(e_2)=0,\quad
	\delta_\prec(e_1)=ae_2\ot e_2,\quad
	\delta_\succ(e_1)=be_2\ot e_2,\quad
	\delta_\bullet(e_1)=ce_2\ot e_2.
	$$

	\item  For all $a,b\in\RR$ with $a\neq0$,
	$$\delta_\prec=\delta_\succ=0,\quad
	\delta_\bullet(e_1)
	=a(e_1\ot e_2+e_2\ot e_1)+be_2\ot e_2,
	\quad\delta_\bullet(e_2)=ae_2\ot e_2.
$$

	\item For all $a,b\in\RR$ with $a\neq0$,
	\[
	\begin{aligned}
		\delta_\prec(e_1)&=ae_1\ot e_2,\quad
		\delta_\succ(e_1)=ae_2\ot e_1,\\
		\delta_\prec(e_2)&=\delta_\succ(e_2)
		=-\delta_\bullet(e_2)=ae_2\ot e_2,\\
		\delta_\bullet(e_1)
		&=-a(e_1\ot e_2+e_2\ot e_1)+be_2\ot e_2.
	\end{aligned}
	\]
\end{enumerate}
\end{ex}

%%%%%%%%%%%%
%%%%%%%%%%%%

\section{Triangular NS-bialgebras, the NS-Yang-Baxter equation,
	and $\mathcal O$-operators}
\mlabel{sec:quasi-triangular-ns}

In this section, we introduce the NS-Yang-Baxter equation
(NS-YBE) and show that its symmetric solutions induce
triangular NS-bialgebras. We then establish the relationship
between symmetric solutions of the NS-YBE and $\mathcal O$-operators
of NS-algebras. Finally, we introduce pre-NS-algebras, study
their correspondence with $\mathcal O$-operators, and construct
nondegenerate symmetric solutions of the NS-YBE from pre-NS-algebras.

\subsection{$\mathcal{O}$-operators of NS-algebras and the NS-Yang-Baxter equation}
We now give the notion of an $\mathcal{O}$-operator of NS-algebras.

\begin{defi}
	Let $(A,\prec,\succ,\bullet)$ be an
	NS-algebra and let
	$r=\sum_i u_i\ot v_i\in A\ot A$.
	For $\circ\in\{\prec,\succ,\bullet,\star\}$,
	we use the following tensor notation:
	{\small 
	\begin{align*}
		r_{12}\circ r_{13}
		:=\sum_{i,j}
		(u_i\circ u_j)\ot v_i\ot v_j,\,
		r_{23}\circ r_{12}
		:=\sum_{i,j}
		u_j\ot(u_i\circ v_j)\ot v_i,\,
		r_{13}\circ r_{23}
		:=\sum_{i,j}
		u_i\ot u_j\ot(v_i\circ v_j).
	\end{align*}
}
	The following system of equations
	is called the {\bf NS-Yang-Baxter
		equation} (or simply the {\bf NS-YBE}):
	\begin{align}
		r_{12}\star r_{13}-r_{23}\succ r_{12}-r_{13}\prec r_{23}&=0,\label{eq:NSYBE1}\\
		r_{12}\bullet r_{13}+r_{23}\bullet r_{12}+	r_{13}\bullet r_{23}&=0.\label{eq:NSYBE2}
	\end{align}
\end{defi}
Equivalently, the NS-Yang-Baxter equation can be written as
\begin{align}
	&\sum_{i,j}
	\bigl(
	(u_i\star u_j)\ot v_i\ot v_j
	-u_j\ot(u_i\succ v_j)\ot v_i
	-u_i\ot u_j\ot(v_i\prec v_j)
	\bigr)=0,
	\mlabel{eq:eybe1}\\
	&\sum_{i,j}
	\bigl(
	(u_i\bullet u_j)\ot v_i\ot v_j
	+u_j\ot(u_i\bullet v_j)\ot v_i
	+u_i\ot u_j\ot(v_i\bullet v_j)
	\bigr)=0.
	\mlabel{eq:eybe2}
\end{align}

\begin{defi}
	Let $(A,\prec,\succ,\bullet)$ be an NS-algebra,
and let
$(V,l_{\prec},r_{\prec},l_{\succ},r_{\succ},
l_{\bullet},r_{\bullet})$
be a representation of $A$. A linear map
	$T: V\to A$ is called an {\bf $\mathcal{O}$-operator} of  $(A,\prec,\succ,\bullet)$  associated  to $(V,l_{\prec},r_{\prec},l_{\succ},r_{\succ},
		l_{\bullet},r_{\bullet})$, if $T$ satisfies 
	\begin{equation}
		T(u)\circ T(v)=T(\ell_\circ(T(u))v+r_\circ(T(v))u),\quad \circ\in \{\prec,\succ,\bullet\},\quad  u,v\in V.
		\mlabel{eq:ope}
	\end{equation}
\end{defi}

\begin{lem}\mlabel{lem:insabt}
Let $(A,\prec,\succ,\bullet)$ be an NS-algebra,
and let $(V,l_{\prec},r_{\prec},l_{\succ},r_{\succ},
l_{\bullet},r_{\bullet})$
be a representation of $(A,\prec,\succ,\bullet)$. Let $T: V\to A$ be an $\mathcal{O}$-operator of  $(A,\prec,\succ,\bullet)$ and $(A\oplus V,\prec_{\ltimes},\succ_{\ltimes},
\bullet_{\ltimes})$ be a semidirect product NS-algebra associated with $(V,l_{\prec},r_{\prec},l_{\succ},r_{\succ},
l_{\bullet},r_{\bullet})$. Define the graph of $T$ by
$\Gamma_T=\{(Tu,u)\mid u\in V\}$. Then $\Gamma_T$ is an NS-subalgebra of $(A\oplus V,\prec_{\ltimes},\succ_{\ltimes},
\bullet_{\ltimes})$, and $(V,\prec_T,\succ_T,\bullet_T)$ is an NS-algebra isomorphic to $\Gamma_T$, where 
$$u\circ_T v:=\ell_\circ(T(u))v+r_\circ(T(v))u,\quad \circ\in \{\prec,\succ,\bullet\},\quad  u,v\in V.$$
\end{lem}
\begin{proof}
	By Eq.~\eqref{eq:spnsa}, for all $u,v\in V$ and for each $\circ\in \{\prec,\succ,\bullet\}$, we have 
	\begin{align}
		(Tu,u)\circ_{\ltimes}(Tv,v)
		&=(Tu\circ Tv,l_{\circ}(Tu)v+r_{\circ}(Tv)u)\notag\\
		&=(T(\ell_\circ(T(u))v+r_\circ(T(v))u),l_{\circ}(Tu)v+r_{\circ}(Tv)u)\notag\\
		&=(T(u\circ_T v),u\circ_T v)\in \Gamma_T.\mlabel{eq:Tuv}
		\end{align}
	Thus, $\Gamma_T$ is an NS-subalgebra of $(A\oplus V,\prec_{\ltimes},\succ_{\ltimes},
	\bullet_{\ltimes})$.

Define $f:V\to\Gamma_T$ by
$f(u)=(T(u),u)$.
Then $f$ is a linear isomorphism
with inverse
$f^{-1}(T(u),u)=u$.
Moreover, for each
$\circ\in\{\prec,\succ,\bullet\}$,
$$f(u\circ_T v)=(T(u\circ_T v),u\circ_T v)\stackrel{\meqref{eq:Tuv}}{=}(Tu,u)\circ_{\ltimes}(Tv,v)=f(u)\circ_{\ltimes}f(v).$$
Hence, by Definition~\mref{defi:NSAI}, $f$ is an NS-algebra isomorphism,
and $(V,\prec_T,\succ_T,\bullet_T)$
is an NS-algebra.
\end{proof}

\begin{pro}\mlabel{pro:6.5}
	Let $(A,\prec,\succ,\bullet)$ be an NS-algebra and
	let $r\in A\ot A$ be symmetric.
	Define $r^\sharp:A^*\to A$ by
\begin{equation}
	\langle\xi,r^\sharp(\eta)\rangle
	=\langle\eta\ot\xi,r\rangle.\mlabel{eq:rsharp}
	\end{equation}
		Then the following conditions
	are equivalent:
	\begin{enumerate}
		\item $r$ is a solution of
	 the NS-YBE.
		\item $r^\sharp$ is an
		$\mathcal O$-operator of $(A,\prec,\succ,\bullet)$
		associated with 
		$
		(A^*,-R_\succ^*,L_\star^*,
		R_\star^*,-L_\prec^*,
		-R_\bullet^*,-L_\bullet^*).
		$
	\end{enumerate}
\end{pro}
\begin{proof}
	Write
$
	r=\sum_i u_i\ot v_i,
$
	and put $T:=r^\sharp:A^*\to A$.
	Since $r$ is symmetric, for any
	$\xi\in A^*$, we have
	\begin{equation}
		T(\xi)
		\stackrel{\meqref{eq:rsharp}}{=}\sum_i\langle\xi,u_i\rangle v_i
		=\sum_i\langle\xi,v_i\rangle u_i.
		\label{eq:NSYBE-sharp}
	\end{equation}
By the definition of an
	$\mathcal O$-operator associated
	with the coadjoint representation,
	$T$ is an $\mathcal O$-operator
	if and only if, for all
	$\xi,\eta\in A^*$,
	\begin{align}
		T(\xi)\prec T(\eta)
		&=T\bigl(
		-R_\succ^*(T(\xi))\eta
		+L_\star^*(T(\eta))\xi
		\bigr),
		\label{eq:NSYBE-op1}\\
		T(\xi)\succ T(\eta)
		&=T\bigl(
		R_\star^*(T(\xi))\eta
		-L_\prec^*(T(\eta))\xi
		\bigr),
		\label{eq:NSYBE-op2}\\
		T(\xi)\bullet T(\eta)
		&=T\bigl(
		-R_\bullet^*(T(\xi))\eta
		-L_\bullet^*(T(\eta))\xi
		\bigr).
		\label{eq:NSYBE-op3}
	\end{align}
For all $\xi,\eta,\zeta\in A^*$, we have
\begin{align*}
\langle\zeta,T(\xi)\prec T(\eta)\rangle
&\stackrel{\meqref{eq:NSYBE-sharp}}{=}\langle\zeta,(\sum_i\langle\xi,u_i\rangle v_i)\prec(\sum_j\langle\eta,u_j\rangle v_j) \rangle\\
&=\langle\zeta,\sum_{i,j}\langle\xi,u_i \rangle\langle\eta,u_j \rangle(v_i\prec v_j) \rangle\\
&=\sum_{i,j}\langle\xi,u_i \rangle\langle\eta,u_j \rangle\langle\zeta,v_i\prec v_j \rangle\\
&=\langle\xi\ot\eta\ot\zeta, \sum_{i,j}u_i\ot u_j\ot (v_i\prec v_j)\rangle\\
&=\langle\xi\ot\eta\ot\zeta,r_{13}\prec r_{23} \rangle,
\end{align*}and\begin{align*}
&\quad\langle \zeta,T(-R_\succ^*(T(\xi))\eta+L_\star^*(T(\eta))\xi) \rangle\\
&=\langle\zeta,\sum_i\langle-R_\succ^*(T(\xi))\eta,u_i\rangle v_i \rangle+\langle\zeta,\sum_i\langle L_\star^*(T(\eta))\xi,u_i\rangle v_i \rangle\\
&=-\sum_i\langle R_\succ^*(T(\xi))\eta,u_i\rangle\langle\zeta,v_i \rangle+\sum_i\langle L_\star^*(T(\eta))\xi,u_i \rangle\langle\zeta,v_i \rangle\\
&=-\sum_i\langle\eta,u_i\succ(T(\xi)) \rangle\langle\zeta,v_i \rangle+\sum_i\langle\xi,(T(\eta))\star u_i \rangle\langle\zeta,v_i \rangle\\
&=-\sum_i\langle\eta,u_i\succ(\sum_j\langle\xi,u_j\rangle v_j ) \rangle\langle\zeta,v_i \rangle+\sum_i\langle\xi,( \sum_j\langle\eta,v_j\rangle u_j)\star u_i \rangle\langle\zeta,v_i \rangle\\
&=-\sum_{i,j}\langle\xi,u_j \rangle\langle\eta,u_i\succ v_j \rangle\langle\zeta,v_i \rangle+\sum_{i,j}\langle\eta,v_j \rangle\langle\xi,u_j \star u_i \rangle\langle\zeta,v_i \rangle\\
&=-\langle\xi\ot\eta\ot\zeta,\sum_{i,j}(u_j\ot(u_i\succ v_j)\ot v_i) \rangle+\langle\xi\ot\eta\ot\zeta,\sum_{i,j}((u_i\star u_j)\ot v_i\ot v_j) \rangle\\
&=\langle\xi\ot\eta\ot\zeta, r_{12}\star r_{13}-r_{23}\succ r_{12}\rangle,
\end{align*}
which shows that Eq.~\eqref{eq:NSYBE1} holds if and only if Eq.~\eqref{eq:NSYBE-op1} holds.  Similarly, we have
\begin{eqnarray*}
	\text{Eq}.~\meqref{eq:NSYBE1}\Leftrightarrow \text{Eq}.~\meqref{eq:NSYBE-op2},\quad
	\text{Eq}.~\meqref{eq:NSYBE2}\Leftrightarrow \text{Eq}.~\meqref{eq:NSYBE-op3}.
\end{eqnarray*}
Therefore, $r$ is a solution of
Eqs.~\eqref{eq:NSYBE1}--
\eqref{eq:NSYBE2} if and only if $r^\sharp$ is an
$\mathcal O$-operator of $A$
associated with 
$
(A^*,-R_\succ^*,L_\star^*,
R_\star^*,-L_\prec^*,
-R_\bullet^*,-L_\bullet^*).
$
\end{proof}

\begin{pro}\mlabel{pro:6.2}
	Let $(A,\prec,\succ,\bullet)$ be an NS-algebra and let
	$r=\sum_i u_i\ot v_i\in A\ot A$
	be a symmetric solution of the NS-YBE.
	Then
	$
	(A,\prec,\succ,\bullet,
	\Delta_{\prec},\Delta_{\succ},
	\Delta_{\bullet})
$
	is an NS-bialgebra, where
	\begin{align}
		\Delta_{\prec}(x)
		&=\bigl(R_\star(x)\ot\id-\id\ot L_\succ(x)\bigr)r,\label{eq:doins01}\\
		\Delta_{\succ}(x)
		&=\bigl(\id\ot L_\star(x)-R_\prec(x)\ot\id\bigr)r,\label{eq:doins02}\\
		\Delta_{\bullet}(x)
		&=-\bigl(R_\bullet(x)\ot\id+\id\ot L_\bullet(x)\bigr)r,\label{eq:doins03}
	\end{align}
	for all $x\in A$.
	In this case,
	$(A,\prec,\succ,\bullet,
	\Delta_{\prec},\Delta_{\succ},
	\Delta_{\bullet})$
	is called a {\bf triangular NS-bialgebra}.
\end{pro}
\begin{proof}
		Let $T:=r^\sharp:A^*\to A$.
	Since $r$ is a symmetric solution
	of the NS-YBE, by Proposition~\mref{pro:6.5}, $T$ is an
	$\mathcal O$-operator associated
	with $(A^*,-R_\succ^*,L_\star^*,
	R_\star^*,-L_\prec^*,
	-R_\bullet^*,-L_\bullet^*).$ 
	By Lemma~\mref{lem:insabt}, $(A^*,\prec_r,\succ_r,\bullet_r)$ is an NS-algebra, where 
	\begin{align*}
\xi\prec_r \eta&:=-R_\succ^*(T(\xi))\eta+L_\star^*(T(\eta))\xi,  \\
\xi\succ_r \eta&:=R_\star^*(T(\xi))\eta-L_\prec^*(T(\eta))\xi,\\
\xi\bullet_r \eta&:=-R_\bullet^*(T(\xi))\eta-L_\bullet^*(T(\eta))\xi,\quad\xi,\eta\in A^*.
	\end{align*}
	For each $\circ\in \{\prec,\succ,\bullet\}$, let $\Delta_{\circ}:A\to A\otimes A$
	be the dual linear map of $\circ_r$,
	as defined by Eq.~\meqref{eq:lds1}. By direct computation, $\Delta_{\prec},\Delta_{\succ},
	\Delta_{\bullet}$ are given by Eqs.~\eqref{eq:doins01}--
	\eqref{eq:doins03}, respectively.
Then by Proposition~\mref{pro:dual-ns-coalgebra}, $(A,\Delta_{\prec},\Delta_{\succ},
	\Delta_{\bullet})$ is an NS-coalgebra.  Let $(A,\star)$ be the associated associative algebra of
	$(A,\prec,\succ,\bullet)$. Since $r$ is symmetric, for all $x,y\in A$, we have 	
	\begin{eqnarray*}
&\quad&(L_\star(x)\ot\id)\Delta_\prec(y)+\tau(L_\succ(y)\ot\id)\Delta_\succ(x)\\
&=&(L_\star(x)\ot\id)(R_\star(y)\ot\id-\id\ot L_\succ(y))r+\tau(L_\succ(y)\ot\id)(\id\ot L_\star(x)-R_\prec(x)\ot\id )r\\
&=&(L_\star(x)R_\star(y)\ot\id)r-(L_\star(x)\ot L_\succ(y))r+\tau(L_\succ(y)\ot L_\star(x))r-\tau(L_\succ(y)R_\prec(x)\ot\id)r\\
&\stackrel{(\mref{eq:ronsa05})}=&(L_\star(x)R_\star(y)\ot\id)r-(L_\star(x)\ot L_\succ(y))r+\tau(L_\succ(y)\ot L_\star(x))r-\tau(R_\prec(x)L_\succ(y)\ot\id)r\\
&\stackrel{\text{Prop}.\mref{pro:roaa}}{=}&(R_\star(y)L_\star(x)\ot\id)r-(L_\star(x)\ot L_\succ(y))r+( L_\star(x)\ot L_\succ(y))\tau(r)- (\id\ot R_\prec(x)L_\succ(y))\tau(r)\\
&=&(R_\star(y)L_\star(x)\ot\id)r-(L_\star(x)\ot L_\succ(y))r+( L_\star(x)\ot L_\succ(y))r- (\id\ot R_\prec(x)L_\succ(y))r\\
&=&(R_\star(y)L_\star(x)\ot\id)r- (\id\ot R_\prec(x)L_\succ(y))r,
	\end{eqnarray*}
	and
	\begin{eqnarray*}
&\quad&\tau(\id\ot R_\star(y))\Delta_\succ(x)+(\id\ot R_\prec(x))\Delta_\prec(y)\\
&=&\tau(\id\ot R_\star(y))(\id\ot L_\star(x)-R_\prec(x)\ot\id )r+(\id\ot R_\prec(x))(R_\star(y)\ot\id-\id\ot L_\succ(y))r\\
&=&\tau(\id\ot R_\star(y)L_\star(x))r-\tau(R_\prec(x)\ot R_\star(y))r+(R_\star(y)\ot R_\prec(x))r-(\id\ot R_\prec(x)L_\succ(y))r\\
&=&( R_\star(y)L_\star(x)\ot\id)\tau(r)-(R_\star(y)\ot R_\prec(x))\tau(r)+(R_\star(y)\ot R_\prec(x))r-(\id\ot R_\prec(x)L_\succ(y))r\\
&=&( R_\star(y)L_\star(x)\ot\id)r-(R_\star(y)\ot R_\prec(x))r+(R_\star(y)\ot R_\prec(x))r-(\id\ot R_\prec(x)L_\succ(y))r\\
&=&( R_\star(y)L_\star(x)\ot\id)r-(\id\ot R_\prec(x)L_\succ(y))r,
\end{eqnarray*}
which shows that Eq.~(\mref{eq:NSb06}) holds. Similarly, Eqs.~(\mref{eq:NSb04})--(\mref{eq:NSb05}) hold. Therefore, by Definition~\mref{def:ns-bialgebra}, $(A,\prec,\succ,\bullet,
\Delta_{\prec},\Delta_{\succ},
\Delta_{\bullet})$
is an NS-bialgebra.
\end{proof}
We next prove that $\mathcal{O}$-operators provide numerous solutions of the NS-YBE in the semidirect
product NS-algebra, which, via Proposition~\mref{pro:6.2}, induce NS-bialgebras. 

\begin{thm}\mlabel{thm:semi-dybe}
	Let $(A,\prec_A,\succ_A,\bullet_A)$ be an NS-algebra,
	and let
	$(V,l_{\prec},r_{\prec},l_{\succ},r_{\succ},
	l_{\bullet},r_{\bullet})$
	be a representation of $A$.
	Let $T:V\to A$ be a linear map, regarded as
	an element of $A\otimes V^*$. Then
$
	r=T+\tau(T)
$
	is a symmetric solution of the NS-YBE in
$
	A\ltimes^{-r_{\succ}^*,r_{\star}^*,-r_{\bullet}^*}
	_{l_{\star}^*,-l_{\prec}^*,-l_{\bullet}^*}V^*
$
	if and only if $T$ is an $\mathcal O$-operator
	associated with
	$(V,l_{\prec},r_{\prec},l_{\succ},r_{\succ},
	l_{\bullet},r_{\bullet})$.
\end{thm}

\begin{proof}
	Let $\{e_1,\ldots,e_n\}$ be a basis of $V$
	and $\{e_1^*,\ldots,e_n^*\}$ the corresponding
	dual basis of $V^*$. Then
$
	T=\sum_{i=1}^n T(e_i)\otimes e_i^*.
$
	Consequently,
$
	r=T+\tau(T)
	=\sum_{i=1}^n
	\big(e_i^*\otimes T(e_i)
	+T(e_i)\otimes e_i^*\big).
$
By Eq.~\meqref{eq:spnsa}, the binary operations $\prec, \succ,\bullet$ on $
A\ltimes^{-r_{\succ}^*,r_{\star}^*,-r_{\bullet}^*}
_{l_{\star}^*,-l_{\prec}^*,-l_{\bullet}^*}V^*
$ are given as follows:
	\begin{align*}
		(x+\xi)\prec (y+\eta)
		&=x\prec_A y-r_\succ^*(x)\eta
		+l_\star^*(y)\xi,\\
		(x+\xi)\succ (y+\eta)
		&=x\succ_A y+r_\star^*(x)\eta
		-l_\prec^*(y)\xi,\\
		(x+\xi)\bullet (y+\eta)
		&=x\bullet_A y-r_\bullet^*(x)\eta
		-l_\bullet^*(y)\xi.
	\end{align*}
	In particular,
	\[
	\begin{aligned}
		e_i^*\prec  T(e_j)&=l_\star^*(T(e_j))e_i^*,\quad
		T(e_i)\prec  e_j^*=-r_\succ^*(T(e_i))e_j^*,\quad
		e_i^*\succ T(e_j)=-l_\prec^*(T(e_j))e_i^*,\\
		T(e_i)\succ e_j^*&=r_\star^*(T(e_i))e_j^*,\quad
		e_i^*\bullet T(e_j)=-l_\bullet^*(T(e_j))e_i^*,\quad
		T(e_i)\bullet e_j^*=-r_\bullet^*(T(e_i))e_j^*.
	\end{aligned}
	\]
	
	We first consider Eq.~\meqref{eq:NSYBE1}, that is,
	\[
	r_{12}\star r_{13}
	-r_{23}\succ r_{12}
	-r_{13}\prec r_{23}=0.
	\]
By direct computation, we have
	\begin{align*}
		r_{12}\star r_{13}
		={}&\sum_{i,j}
		e_i^*\star e_j^*\otimes T(e_i)\otimes T(e_j)+\sum_{i,j}
		e_i^*\star T(e_j)\otimes T(e_i)\otimes e_j^*\\
		&+\sum_{i,j}
		T(e_i)\star e_j^*\otimes e_i^*\otimes T(e_j)+\sum_{i,j}
		T(e_i)\star T(e_j)\otimes e_i^*\otimes e_j^*\\
		={}&\sum_{i,j}
		l_\succ^*(T(e_j))e_i^*
		\otimes T(e_i)\otimes e_j^*+\sum_{i,j}
		r_\prec^*(T(e_i))e_j^*
		\otimes e_i^*\otimes T(e_j)\\
		&+\sum_{i,j}
		T(e_i)\star T(e_j)
		\otimes e_i^*\otimes e_j^*.
	\end{align*}
	
	Similarly,
	\begin{align*}
		r_{23}\succ r_{12}
		={}&\sum_{i,j}
		e_j^*\otimes
		(e_i^*\succ T(e_j))\otimes T(e_i)+\sum_{i,j}
		e_j^*\otimes
		(T(e_i)\succ T(e_j))\otimes e_i^*\\
		&+\sum_{i,j}
		T(e_j)\otimes
		(T(e_i)\succ e_j^*)\otimes e_i^*\\
		={}&-\sum_{i,j}
		e_j^*\otimes
		l_\prec^*(T(e_j))e_i^*\otimes T(e_i)+\sum_{i,j}
		e_j^*\otimes
		(T(e_i)\succ T(e_j))\otimes e_i^*\\
		&+\sum_{i,j}
		T(e_j)\otimes
		r_\star^*(T(e_i))e_j^*\otimes e_i^*.
	\end{align*}
	
	Moreover,
	\begin{align*}
		r_{13}\prec r_{23}
		={}&\sum_{i,j}
		e_i^*\otimes e_j^*
		\otimes(T(e_i)\prec T(e_j))+\sum_{i,j}
		e_i^*\otimes T(e_j)
		\otimes(T(e_i)\prec e_j^*)\\
		&+\sum_{i,j}
		T(e_i)\otimes e_j^*
		\otimes(e_i^*\prec T(e_j))\\
		={}&\sum_{i,j}
		e_i^*\otimes e_j^*
		\otimes(T(e_i)\prec T(e_j))-\sum_{i,j}
		e_i^*\otimes T(e_j)
		\otimes R_\succ^*(T(e_i))e_j^*\\
		&+\sum_{i,j}
		T(e_i)\otimes e_j^*
		\otimes L_\star^*(T(e_j))e_i^*.
	\end{align*}
Substituting the above expressions into the first
	NS-YBE and comparing the components in
	\[
	V^*\otimes V^*\otimes A,\qquad
	V^*\otimes A\otimes V^*,\qquad
	A\otimes V^*\otimes V^*,
	\]
	respectively, we obtain the following three
	equations. The $V^*\otimes V^*\otimes A$ component is
	\begin{align}
		\sum_{i,j}\big(
		r_\prec^*(T(e_i))e_j^*
		\otimes e_i^*\otimes T(e_j)
		+e_j^*\otimes
		l_\prec^*(T(e_j))e_i^*\otimes T(e_i)
		-e_i^*\otimes e_j^*
		\otimes(T(e_i)\prec T(e_j))
		\big)=0.
		\mlabel{eq:ns-t1}
	\end{align}
	By definition,
	\[
	\begin{aligned}
		\sum_{i,j}
		r_\prec^*(T(e_i))e_j^*
		\otimes e_i^*\otimes T(e_j)
		&=\sum_{i,j,k}
		\langle r_\prec^*(T(e_i))e_j^*,e_k\rangle
		e_k^*\otimes e_i^*\otimes T(e_j)\\
		&=\sum_{i,j,k}
		\langle e_j^*,r_\prec(T(e_i))e_k\rangle
		e_k^*\otimes e_i^*\otimes T(e_j)\\
		&=\sum_{i,k}
		e_k^*\otimes e_i^*\otimes
		T(r_\prec(T(e_i))e_k)\\
		&=\sum_{i,j}
		e_i^*\otimes e_j^*\otimes
		T(r_\prec(T(e_j))e_i).
	\end{aligned}
	\]
	Similarly,
	\[
	\begin{aligned}
		\sum_{i,j}
		e_j^*\otimes
		l_\prec^*(T(e_j))e_i^*\otimes T(e_i)
		&=\sum_{i,j,k}
		e_j^*\otimes
		\langle e_i^*,l_\prec(T(e_j))e_k\rangle
		e_k^*\otimes T(e_i)\\
		&=\sum_{i,j}
		e_i^*\otimes e_j^*\otimes
		T(l_\prec(T(e_i))e_j).
	\end{aligned}
	\]
	Thus, Eq.~\meqref{eq:ns-t1} is equivalent to
	\[
	\begin{aligned}
		T(e_i)\prec T(e_j)
		=T\big(l_\prec(T(e_i))e_j
		+r_\prec(T(e_j))e_i\big).
	\end{aligned}
	\]
	The $V^*\otimes A\otimes V^*$ component is
	\begin{align}
		\sum_{i,j}\big(
		l_\succ^*(T(e_j))e_i^*
		\otimes T(e_i)\otimes e_j^*
		-e_j^*\otimes
		(T(e_i)\succ T(e_j))\otimes e_i^*
		+e_i^*\otimes T(e_j)
		\otimes R_\succ^*(T(e_i))e_j^*
		\big)=0.
		\mlabel{eq:ns-t2}
	\end{align}
	Applying the dual pairing to the first and
	third tensor factors, we obtain
	\[
	\begin{aligned}
		\sum_{i,j}
		l_\succ^*(T(e_j))e_i^*
		\otimes T(e_i)\otimes e_j^*
		&=\sum_{i,j,k}
		\langle e_i^*,l_\succ(T(e_j))e_k\rangle
		e_k^*\otimes T(e_i)\otimes e_j^*\\
		&=\sum_{j,k}
		e_k^*\otimes T(l_\succ(T(e_j))e_k)
		\otimes e_j^*.
	\end{aligned}
	\]
	Likewise,
	\[
	\begin{aligned}
		\sum_{i,j}
		e_i^*\otimes T(e_j)
		\otimes R_\succ^*(T(e_i))e_j^*
		&=\sum_{i,j,k}
		e_i^*\otimes T(e_j)\otimes
		\langle e_j^*,r_\succ(T(e_i))e_k\rangle e_k^*\\
		&=\sum_{i,k}
		e_i^*\otimes T(r_\succ(T(e_i))e_k)
		\otimes e_k^*.
	\end{aligned}
	\]
	Consequently, Eq.~\meqref{eq:ns-t2} is equivalent to
	\[
	\begin{aligned}
		T(e_i)\succ T(e_j)
		=T\big(l_\succ(T(e_i))e_j
		+r_\succ(T(e_j))e_i\big).
	\end{aligned}
	\]
	Finally, the $A\otimes V^*\otimes V^*$ component is
	\begin{align}
		\sum_{i,j}\big(
		T(e_i)\star T(e_j)\otimes e_i^*\otimes e_j^*-T(e_j)\otimes
		r_\star^*(T(e_i))e_j^*\otimes e_i^*-T(e_i)\otimes e_j^*
		\otimes l_\star^*(T(e_j))e_i^*
		\big)=0.
		\mlabel{eq:ns-t3}
	\end{align}
	Likewise, we obtain	\[
	\begin{aligned}
		\sum_{i,j}
		T(e_j)\otimes
		r_\star^*(T(e_i))e_j^*\otimes e_i^*
		=\sum_{i,j}
		T(r_\star(T(e_j))e_i)\otimes e_i^*\otimes e_j^*.
	\end{aligned}
	\]
	Similarly,
	\[
	\begin{aligned}
		\sum_{i,j}
		T(e_i)\otimes e_j^*
		\otimes l_\star^*(T(e_j))e_i^*
		&=\sum_{i,j}
		T(l_\star(T(e_i))e_j)
		\otimes e_i^*\otimes e_j^*.
	\end{aligned}
	\]
	Therefore, Eq.~\meqref{eq:ns-t3} is equivalent to
	\[
	T(e_i)\star T(e_j)
	=T\big(
	l_\star(T(e_i))e_j+
	r_\star(T(e_j))e_i
	\big).
	\]
	
	Since
	$
	\star=\prec+\succ+\bullet,
	$
	and the identities for $\prec$ and $\succ$
	have already been established, we conclude that
	\[
	T(e_i)\bullet T(e_j)
	=T\big(
	l_\bullet(T(e_i))e_j+
	r_\bullet(T(e_j))e_i
	\big).
	\]
	
	Thus, Eq.~\meqref{eq:NSYBE1} yields Eq.~\meqref{eq:ope}.
	
	We now consider Eq.~\meqref{eq:NSYBE2}, that is,
	\[
	r_{12}\bullet r_{13}
	+r_{23}\bullet r_{12}
	+r_{13}\bullet r_{23}=0.
	\]
	
	Direct computation gives
	{\small 
	\begin{align*}
		r_{12}\bullet r_{13}
		={}-\sum_{i,j}
		l_\bullet^*(T(e_j))e_i^*
		\otimes T(e_i)\otimes e_j^*
		-\sum_{i,j}
		r_\bullet^*(T(e_i))e_j^*
		\otimes e_i^*\otimes T(e_j)
		+\sum_{i,j}
		T(e_i)\bullet T(e_j)
		\otimes e_i^*\otimes e_j^*,
		\\
		r_{23}\bullet r_{12}
		={}-\sum_{i,j}
		e_j^*\otimes
		l_\bullet^*(T(e_j))e_i^*
		\otimes T(e_i)+\sum_{i,j}
		e_j^*\otimes
		(T(e_i)\bullet T(e_j))\otimes e_i^*
		-\sum_{i,j}
		T(e_j)\otimes
		r_\bullet^*(T(e_i))e_j^*
		\otimes e_i^*,
		\\
		r_{13}\bullet r_{23}
		={}\sum_{i,j}
		e_i^*\otimes e_j^*
		\otimes(T(e_i)\bullet T(e_j))
		-\sum_{i,j}
		e_i^*\otimes T(e_j)
		\otimes r_\bullet^*(T(e_i))e_j^*
		-\sum_{i,j}
		T(e_i)\otimes e_j^*
		\otimes l_\bullet^*(T(e_j))e_i^*.
	\end{align*}
}	

A similar decomposition of
Eq.~\meqref{eq:NSYBE2} shows that
its three nonzero tensor components
are equivalent to the
$\mathcal O$-operator identity
for $\bullet$, which already follows
from Eq.~\meqref{eq:NSYBE1}.

	Conversely, suppose that $T$ is an
	$\mathcal O$-operator associated with the
	given representation. Then the three
	$\mathcal O$-operator identities hold.
	By the above computations, all the
	components of the two NS-YBE equations
	vanish. Hence $r=T+\tau(T)$ satisfies
	the NS-YBE in the semidirect product NS-algebra.
	Therefore, $r=T+\tau(T)$ is a symmetric
	solution of the NS-YBE if and only if
	$T$ is an $\mathcal O$-operator.
\end{proof}

\subsection{Pre-NS-algebras and $\mathcal{O}$-operators}

\begin{defi}\mlabel{def:pre-ns}
	A {\bf pre-NS-algebra} is a vector
	space $V$ equipped with six bilinear
	operations
	\[
	\lhd_{\prec},\ \rhd_{\prec},\
	\lhd_{\succ},\ \rhd_{\succ},\
	\lhd_{\bullet},\ \rhd_{\bullet}
	\]
	such that $L_{\lhd_\prec}, R_{\rhd_\prec}, L_{\lhd_\succ}, R_{\rhd_\succ},L_{\lhd_\bullet}, R_{\rhd_\bullet}$ satisfy Eqs.~\meqref{eq:ronsa01}--\meqref{eq:ronsa12}. Equivalently, 
 for all $u,v, w\in V$,
\begin{align}
	(u\prec_V v)\lhd_\prec w
	&=u\lhd_\prec(v\lhd_\star w),
	\label{eq:pnsa01}\\
		(u\lhd_\prec v)\rhd_\prec w
&=u\lhd_\prec(v\rhd_\star w),
\label{eq:pnsa02}\\
(u\rhd_\prec v)\rhd_\prec w
&=u\rhd_\prec(v\star_V w),
\label{eq:pnsa03}\\
(u\succ_V v)\lhd_\prec w
&=u\lhd_\succ(v\lhd_\prec w),
\label{eq:pnsa04}\\
(u\lhd_\succ v)\rhd_\prec w
&=u\lhd_\succ(v\rhd_\prec w),
\label{eq:pnsa05}\\
u\rhd_\succ(v\prec_V w)
&=(u\rhd_\succ v)\rhd_\prec w,
\label{eq:pnsa06}\\
u\rhd_\succ(v\succ_V w)
&=(u\rhd_\star v)\rhd_\succ w,
\label{eq:pnsa07}\\
u\lhd_\succ(v\rhd_\succ w)
&=(u\lhd_\star v)\rhd_\succ w,
\label{eq:pnsa08}\\
u\lhd_\succ(v\lhd_\succ w)
&=(u\star_V v)\lhd_\succ w,
\label{eq:pnsa09}\\
(u\star_V v)\lhd_\bullet w
+(u\bullet_V v)\lhd_\prec w
&=u\lhd_\succ(v\lhd_\bullet w)
+u\lhd_\bullet(v\lhd_\star w),
\label{eq:pnsa10}\\
u\lhd_\succ(v\rhd_\bullet w)
+u\lhd_\bullet(v\rhd_\star w)
&=(u\lhd_\bullet v)\rhd_\prec w
+(u\lhd_\star v)\rhd_\bullet w,
\label{eq:pnsa11}\\
u\rhd_\bullet(v\star_V w)
+u\rhd_\succ(v\bullet_V w)
&=(u\rhd_\bullet v)\rhd_\prec w
+(u\rhd_\star v)\rhd_\bullet w.
\label{eq:pnsa12}
\end{align}
Here 	
for each $\circ\in\{\prec,\succ,\bullet\}$,
define
\begin{align}
	u\circ_V v
	:=u\lhd_{\circ} v+u\rhd_{\circ} v,\label{eq:pnsa}
\end{align} and $\star$ is defined by Eq.~\meqref{eq:ast}.
\end{defi}

\begin{pro}\mlabel{pro:pre-ns-ooperator}
	Let $(A,\prec,\succ,\bullet)$
	be an NS-algebra, and let $(V,l_\prec,r_\prec,l_\succ,r_\succ,
	l_\bullet,r_\bullet)$
	be a representation of $A$.
	If $T:V\to A$ is an
	$\mathcal O$-operator associated
	with this representation, then
	the six operations
\begin{align}
	\begin{split}
		u\lhd_\prec v&=l_\prec(Tu)v,\quad
		u\rhd_\prec v=r_\prec(Tv)u,\quad
		u\lhd_\succ v=l_\succ(Tu)v,\\
		u\rhd_\succ v&=r_\succ(Tv)u,\quad
		u\lhd_\bullet v=l_\bullet(Tu)v,\quad
		u\rhd_\bullet v=r_\bullet(Tv)u.
	\end{split}
	\mlabel{eq:ropnsa}
\end{align}

	make $V$ into a pre-NS-algebra.
	Moreover, for each
	$\circ\in\{\prec,\succ,\bullet\}$, the induced operations $u\circ_V v$ by Eq.~\meqref{eq:pnsa}
	make $V$ into an NS-algebra, and
	$T:V\to A$ is an NS-algebra
	homomorphism.
\end{pro}

\begin{proof}
	Since $(A,\prec,\succ,\bullet)$
	is an NS-algebra and $T:V\to A$ is an
	$\mathcal O$-operator associated
	with $(V,l_\prec,r_\prec,l_\succ,r_\succ,
	l_\bullet,r_\bullet)$, by Lemma~\mref{lem:insabt}, the
	induced operations
	$\prec_V,\succ_V,\bullet_V$ make $V$
	into an NS-algebra. Furthermore, for each
	$\circ\in\{\prec,\succ,\bullet\}$, we have 
	\[
	T(u\circ_V v)=T(u\lhd_\circ v+u\rhd_\circ v)=T(l_\circ(T(u))v+r_\circ(T(v))u)\stackrel{(\mref{eq:ope})}=T(u)\circ T(v).
	\]
	Thus, $T$ is a homomorphism
	with respect to the three
	induced operations. For all
	$u,v,w\in V$, by Eq.~(\mref{eq:ropnsa}),
	we have 
	\begin{eqnarray*}
		(u\prec_V v)\lhd_\prec w
		&=&l_\prec(T(u\prec_V v)) w\\
		&=&l_\prec(Tu\prec Tv)w\\
		&\stackrel{(\mref{eq:ronsa01})}=&l_\prec(Tu)l_\star(Tv)w\\
		&=&u\lhd_\prec(v\lhd_\star w),
	\end{eqnarray*}
    which shows that Eq.~(\mref{eq:pnsa01}) holds.
	Similarly, Eqs.~(\mref{eq:pnsa02})--(\mref{eq:pnsa12}) hold. Therefore, $V$ is a pre-NS-algebra.
\end{proof}

\begin{pro}\mlabel{pro:pre-ns-characterization}
	Let $V$ be a vector space equipped with
	binary operations
$
	\lhd_\prec,\ \rhd_\prec,\
	\lhd_\succ,\ \rhd_\succ,\
	\lhd_\bullet,\ \rhd_\bullet.
$
	For each $\circ\in\{\prec,\succ,\bullet\}$,
	define $\circ_V:V\ot V\to V$  by Eq.~(\mref{eq:pnsa}).
	Then the following statements are equivalent:
	\begin{enumerate}
		\item\mlabel{it:pnsa11}
		$(V,\lhd_\prec,\rhd_\prec,
		\lhd_\succ,\rhd_\succ,
		\lhd_\bullet,\rhd_\bullet)$
		is a pre-NS-algebra.
		
		\item\mlabel{it:pnsa12}
$(V,\prec_V,\succ_V,\bullet_V)$ is an NS-algebra, called the
{\bf sub-adjacent NS-algebra} of
$(V,\lhd_\prec,\rhd_\prec,
\lhd_\succ,\rhd_\succ,\\
\lhd_\bullet,\rhd_\bullet)$.
Moreover,
$(V,L_{\lhd_\prec},R_{\rhd_\prec},
L_{\lhd_\succ},R_{\rhd_\succ},
L_{\lhd_\bullet},R_{\rhd_\bullet})$
is a representation of this NS-algebra, and
$\mathrm{id}_V$ is an $\mathcal O$-operator associated
with this representation.
In this case, the given pre-NS-algebra is called a
{\bf compatible pre-NS-algebra} of its sub-adjacent
NS-algebra.
	\end{enumerate}
\end{pro}

\begin{proof}
	($(\mref{it:pnsa11})\Rightarrow(\mref{it:pnsa12})$ )
	Suppose that $(V,\lhd_\prec,\rhd_\prec,
	\lhd_\succ,\rhd_\succ,
	\lhd_\bullet,\rhd_\bullet)$ is a pre-NS-algebra. For all $x,y,z\in V$, we have
	\begin{eqnarray*}
		(x\prec_V y)\prec_V z
		&\stackrel{(\mref{eq:pnsa})}=&(x\prec_V y)\lhd_{\prec} z+(x\lhd_{\prec} y+x\rhd_{\prec} y)\rhd_{\prec} z\\
		&=&L_{\lhd_\prec}(x\prec_V y)z+R_{\rhd_\prec}(z)(L_{\lhd_\prec}(x)y+R_{\rhd_\prec}(y)x)\\
		&\stackrel{(\mref{eq:pnsa01})}=&L_{\lhd_\prec}(x)L_{\lhd_\star}(y)z+R_{\rhd_\prec}(z)L_{\lhd_\prec}(x)y+R_{\rhd_\prec}(z)R_{\rhd_\prec}(y)x\\
		&\stackrel{(\mref{eq:pnsa02})}=&L_{\lhd_\prec}(x)L_{\lhd_\star}(y)z+L_{\lhd_\prec}(x)R_{\rhd_\star}(z)y+R_{\rhd_\prec}(z)
		R_{\rhd_\prec}(y)x\\
		&\stackrel{(\mref{eq:pnsa03})}=&L_{\lhd_\prec}(x)L_{\lhd_\star}(y)z+L_{\lhd_\prec}(x)R_{\rhd_\star}(z)y+R_{\rhd_\prec}(y\star z)x\\
		&=&x\lhd_{\prec}(y\lhd_{\star}z)+x\lhd_{\prec}(y\rhd_{\star}z)+x\rhd_{\prec}(y\star z)\\
		&\stackrel{(\mref{eq:pnsa})}=&x\prec_V(y\star_V z),
	\end{eqnarray*}
which shows that Eq.~(\mref{eq:nsa01}) holds.
Similarly, Eqs.~(\mref{eq:nsa02})--(\mref{eq:nsa04}) hold. Therefore, $(V,\prec_V,\succ_V,\bullet_V)$ is an NS-algebra.
By Definition~\mref{def:pre-ns}, the six linear maps
	$
	L_{\lhd_\prec},R_{\rhd_\prec},L_{\lhd_\succ},R_{\rhd_\succ},
	L_{\lhd_\bullet},R_{\rhd_\bullet}
	$
	satisfy the representation identities
	of an NS-algebra with respect to the
	induced operations $\prec_V,\succ_V,\bullet_V$. Moreover, for each 
	$\circ\in\{\prec,\succ,\bullet\}$ and for all $x,y\in V$,
	we have
	\begin{align*}
		\mathrm{id}_V(x)\circ_V\mathrm{id}_V(y)
	%	&=x\circ_V y\\
		=x\lhd_\circ y+x\rhd_\circ y
		=L_\circ(x)y+R_\circ(y)x
		=\mathrm{id}_V\bigl(
		L_\circ(\mathrm{id}_V(x))y
		+R_\circ(\mathrm{id}_V(y))x
		\bigr).
	\end{align*}
	Thus, $\mathrm{id}_V$ is an
	$\mathcal O$-operator associated
	with $(V,L_{\lhd_\prec},R_{\rhd_\prec},L_{\lhd_\succ},R_{\rhd_\succ},
	L_{\lhd_\bullet},R_{\rhd_\bullet})$.
	
	($(\mref{it:pnsa12})\Rightarrow(\mref{it:pnsa11})$) 
Suppose that
	$(V,\prec_V,\succ_V,\bullet_V)$ is an
	NS-algebra and
	$(V,L_{\lhd_\prec},R_{\rhd_\prec},L_{\lhd_\succ},R_{\rhd_\succ},
	L_{\lhd_\bullet},R_{\rhd_\bullet})$
	is a representation of 	$(V,\prec_V,\succ_V,\bullet_V)$. Then by Definition~\mref{def:pre-ns}, $(V,\lhd_\prec,\rhd_\prec,
	\lhd_\succ,\rhd_\succ,
	\lhd_\bullet,\rhd_\bullet)$
	is a pre-NS-algebra.
\end{proof}

\begin{cor}\mlabel{cor:pre-ns-identity}
Every $\mathcal O$-operator of an NS-algebra
	with respect to a representation induces a pre-NS-algebra.
	Conversely,	every pre-NS-algebra gives rise to an $\mathcal O$-operator.
\end{cor}

\begin{proof}
	This follows from Propositions~\mref{pro:pre-ns-ooperator} and~\mref{pro:pre-ns-characterization}.
\end{proof}
The following is an equivalence between compatible pre-NS-algebra structures on an NS-algebra and invertible $\mathcal{O}$-operators.

\begin{pro}\mlabel{pro:inv-ooperator-pre-ns}
	Let $(A,\prec,\succ,\bullet)$ be an NS-algebra.
	Then $A$ admits a compatible pre-NS-algebra
	structure
$(A,\lhd'_\prec,\rhd'_\prec,
	\lhd'_\succ,\rhd'_\succ,
	\lhd'_\bullet,\rhd'_\bullet)$
	such that
$x\circ y=x\lhd'_\circ y+x\rhd'_\circ y,
\,\circ\in\{\prec,\succ,\bullet\},
$
if and only if there exists an invertible
	$\mathcal O$-operator $T:V\to A$ associated
	with a representation
$
	(V,l_\prec,r_\prec,l_\succ,r_\succ,
	l_\bullet,r_\bullet)
$
	of $(A,\prec,\succ,\bullet)$.
Moreover, if such an operator $T$ exists,
	the compatible pre-NS-algebra operations
	on $A$ are given by
\begin{align*}
	&x \lhd'_{\prec} y := T\left( l_{\prec}(x)T^{-1}(y) \right), \, x \rhd'_\prec y := T\left( r_{\prec}(y) T^{-1}(x) \right),\\
	&x \lhd'_\succ y := T\left( l_{\succ}(x)T^{-1}(y) \right), \, x \rhd'_\succ y := T\left( r_{\succ}(y) T^{-1}(x) \right), \\
	&x \lhd'_\bullet y := T\left( l_{\bullet}(x)T^{-1}(y) \right), \, x \rhd'_\bullet y := T\left( r_{\bullet}(y) T^{-1}(x) \right),\quad x, y \in A. 
\end{align*}
\end{pro}

\begin{proof}
	Suppose first that $A$ admits a compatible
	pre-NS-algebra structure
$(A,\lhd'_\prec,\rhd'_\prec,
	\lhd'_\succ,\rhd'_\succ,
	\lhd'_\bullet,\rhd'_\bullet).$
	By Proposition~\mref{pro:pre-ns-characterization},
$(A,L_{\lhd'_\prec},R_{\rhd'_\prec},
	L_{\lhd'_\succ},R_{\rhd'_\succ},
	L_{\lhd'_\bullet},R_{\rhd'_\bullet})$
	is a representation of
	$(A,\prec,\succ,\bullet)$, and
	$\mathrm{id}_A$ is an $\mathcal O$-operator
	associated with this representation.
	Since $\mathrm{id}_A$ is invertible,
	the necessity follows.
	
	Conversely, suppose that $T:V\to A$
	is an invertible $\mathcal O$-operator
	associated with the representation
$(V,l_\prec,r_\prec,l_\succ,r_\succ,
	l_\bullet,r_\bullet).$
By Proposition~\mref{pro:pre-ns-ooperator},
	the six operations
	\[
	u\lhd_\circ v=l_\circ(Tu)v,
	\qquad
	u\rhd_\circ v=r_\circ(Tv)u,
	\quad
	\circ\in\{\prec,\succ,\bullet\},\,u,v\in V
	\]
	make $V$ into a pre-NS-algebra.
	Since $T$ is invertible, define six
	bilinear operations on $A$ by
	\[
	\begin{aligned}
		x\lhd'_\circ y
		=T\bigl(
		T^{-1}(x)\lhd_\circ T^{-1}(y)
		\bigr),\quad
		x\rhd'_\circ y
		=T\bigl(
		T^{-1}(x)\rhd_\circ T^{-1}(y)
		\bigr),\quad x,y\in A,\,\circ\in\{\prec,\succ,\bullet\}.
	\end{aligned}
	\]
	Equivalently,
	\[
	\begin{aligned}
		x\lhd'_\circ y
		=T\bigl(l_\circ(x)T^{-1}(y)\bigr),\quad
		x\rhd'_\circ y
		=T\bigl(r_\circ(y)T^{-1}(x)\bigr).
	\end{aligned}
	\]
	By construction, $T$ is an isomorphism
	of vector spaces preserving all six
	operations. Hence these operations
	make $A$ into a pre-NS-algebra.
	It remains to verify compatibility with
	the given NS-algebra structure on $A$.
	For each $\circ\in\{\prec,\succ,\bullet\}$
	and all $x,y\in A$, we have
	\[
	\begin{aligned}
		x\lhd'_\circ y+x\rhd'_\circ y
		&=T\bigl(
		l_\circ(x)T^{-1}(y)
		+r_\circ(y)T^{-1}(x)
		\bigr)\\
		&\stackrel{(\mref{eq:ope})}{=}
		T(T^{-1}(x))\circ T(T^{-1}(y))\\
		&=x\circ y.
	\end{aligned}
	\]
	Therefore, the induced NS-algebra
	operations coincide with the original
	operations on $A$. Thus, the pre-NS-algebra
	structure on $A$ is compatible with
	$(A,\prec,\succ,\bullet)$.
\end{proof}
Let $(A,\lhd_\prec,\rhd_\prec,
\lhd_\succ,\rhd_\succ,
\lhd_\bullet,\rhd_\bullet)$ be a pre-NS-algebra.
By Proposition~\mref{pro:pre-ns-characterization},
the induced operations $\prec,\succ,\bullet$
make $A$ into an NS-algebra, and
$
(A,L_{\lhd_\prec},R_{\rhd_\prec},
L_{\lhd_\succ},R_{\rhd_\succ},
L_{\lhd_\bullet},R_{\rhd_\bullet})
$
is a representation of it.
By Propositions~\mref{pro:dul-rep}
and~\mref{pro:spnsa}, the vector space
$A\oplus A^*$ carries a semidirect product
NS-algebra structure
\[
\mathcal D(A):=
A\ltimes^{-R_{\rhd_\succ}^*,
	R_{\rhd_\star}^*,-R_{\rhd_\bullet}^*}
_{L_{\lhd_\star}^*,
	-L_{\lhd_\prec}^*,-L_{\lhd_\bullet}^*}A^*.
\]
The following theorem constructs a
nondegenerate symmetric solution of the
NS-YBE from any 
pre-NS-algebra and determines its
associated bilinear form.

\begin{thm}\mlabel{thm:pre-ns}
	Let $(A,\lhd_\prec,\rhd_\prec,
	\lhd_\succ,\rhd_\succ,
	\lhd_\bullet,\rhd_\bullet)$ be a pre-NS-algebra.
	Let $\{e_1,\ldots,e_n\}$ be a basis
	of $A$, with dual basis
	$\{e_1^*,\ldots,e_n^*\}$.
	Then
$
	r:=\sum_{i=1}^n
	(e_i^*\otimes e_i+e_i\otimes e_i^*)
$
	is a nondegenerate symmetric solution
	of the NS-YBE in $\mathcal D(A)$.
	Moreover, the bilinear form $\omega_r$
	on $\mathcal D(A)$ induced by $r$
	is given by
	\begin{equation}\mlabel{eq:inv-omeg}
		\omega_r(x+a^*,y+b^*)
		=\langle a^*,y\rangle+\langle b^*,x\rangle,
		\quad
		x,y\in A,\quad a^*,b^*\in A^*.
	\end{equation}
\end{thm}

\begin{proof}By Proposition~\mref{pro:pre-ns-characterization},
	$\mathrm{id}_A$ is an $\mathcal O$-operator
	of the sub-adjacent NS-algebra
	associated with its natural representation.
	Hence, by Theorem~\mref{thm:semi-dybe},
	$r$ is a symmetric solution of the NS-YBE
	in $\mathcal D(A)$.
	It remains to prove that $r$ is
	nondegenerate and determine the
	associated bilinear form.
	Let
	$
	r^\sharp:(A\oplus A^*)^*
	\longrightarrow A\oplus A^*
	$
	be the linear map induced by $r$,
	defined by
	\[
	\langle r^\sharp(\alpha),\beta\rangle
	=\langle\alpha\otimes\beta,r\rangle,
	\qquad
	\alpha,\beta\in(A\oplus A^*)^*.
	\]
	Using the canonical identification
	$
	(A\oplus A^*)^*\cong A^*\oplus A,
	$
	we have, for all $a^*\in A^*$ and
	$x\in A$,
	\[
	\begin{aligned}
		r^\sharp(a^*+x)
		&=\sum_{i=1}^n
		\bigl(
		\langle a^*+x,e_i^*\rangle e_i
		+\langle a^*+x,e_i\rangle e_i^*
		\bigr)\\
		&=\sum_{i=1}^n
		\bigl(
		\langle e_i^*,x\rangle e_i
		+\langle a^*,e_i\rangle e_i^*
		\bigr)\\
		&=x+a^*.
	\end{aligned}
	\]
	Thus, $r^\sharp$ is an isomorphism,
	and hence $r$ is nondegenerate.
	The bilinear form $\omega_r$
	associated with $r$ is defined by
	\[
	\omega_r(u,v)
	=\langle(r^\sharp)^{-1}(u),v\rangle,
	\qquad u,v\in A\oplus A^*.
	\]
	Since
	$
	(r^\sharp)^{-1}(x+a^*)=a^*+x,
	$
	we obtain
	\[
	\begin{aligned}
		\omega_r(x+a^*,y+b^*)
		=\langle a^*+x,y+b^*\rangle
		=\langle a^*,y\rangle
		+\langle b^*,x\rangle,\quad x,y\in A,\,a^*,b^*\in A^*.
	\end{aligned}
	\]
	This proves Eq.~\meqref{eq:inv-omeg}
	and completes the proof.
\end{proof}

\section{Applications of NS-bialgebras to dendriform D-bialgebras}
\mlabel{sec:Dbia}

In this section, we explore the relationship between
NS-bialgebras and dendriform D-bialgebras.
We first show that matched pairs of NS-algebras
induce matched pairs of associative algebras
and that NS-coalgebras admit coassociative coalgebras.
We then construct nondegenerate Connes cocycles
from quadratic NS-algebras and establish their
double constructions.
As an application, we show that every NS-bialgebra
gives rise to a dendriform D-bialgebra
on the direct sum of two copies of its underlying space.

\begin{pro}\mlabel{pro:mpoaa}
Let $\bigl( A_1,A_2,
l_{\prec_1},r_{\prec_1},l_{\succ_1},r_{\succ_1},
l_{\bullet_1},r_{\bullet_1},
l_{\prec_2},r_{\prec_2},l_{\succ_2},r_{\succ_2},
l_{\bullet_2},r_{\bullet_2}\bigr)$ be a matched pair of
NS-algebras. Then $\bigl(A_1,
A_2,l_{\star_1},r_{\star_1},l_{\star_2},r_{\star_2} \bigr)$ is a matched pair of
associative algebras, where, for $i=1,2$, the operations $\star_i$ and maps
$l_{\star_i},r_{\star_i}$ are defined by Eqs.~\meqref{eq:ast} and \meqref{eq:lrast}, respectively.
\end{pro}
\begin{proof}By Proposition \mref{pro:roaa}, for $i=1,2$, the algebras $(A_i,\star_i)$
	are associative, and
	$(A_j,l_{\star_i},r_{\star_i})$
	is a representation of $(A_i,\star_i)$ for $i\neq j$.
 By Proposition \mref{pro:matched-ns}, 
$(A_1\oplus A_2,\prec_{\bowtie},\succ_{\bowtie},\bullet_{\bowtie})$ is an NS-algebra, where the three bilinear operations $\prec_{\bowtie},\succ_{\bowtie},
\bullet_{\bowtie}$ on $A_1\oplus A_2$ are given by Eq.~\meqref{eq:bp}. Thus, $(A_1\oplus A_2,\star_{\bowtie})$ is an associative algebra, where the bilinear operation $\star_{\bowtie}$ on $A_1\oplus A_2$ is defined by 
\begin{align*}
	(x,u)\star_{\bowtie}(y,v):&=(x,u)\prec_{\bowtie}(y,v)+(x,u)\succ_{\bowtie}(y,v)+(x,u)\bullet_{\bowtie}(y,v)\\
	&=\bigl(x\star_1 y
	+l_{\star_2}(u)y+r_{\star_2}(v)x,
	u\star_2 v
	+l_{\star_1}(x)v+r_{\star_1}(y)u\bigr), \quad x,y\in A_1,\,u,v\in A_2.
\end{align*}
 Therefore, by \cite[Theorem 2.1.4]{Bai1}, $\bigl(A_1,
 A_2,l_{\star_1},r_{\star_1},l_{\star_2},r_{\star_2} \bigr)$ is a matched pair of associative algebras.
\end{proof}

\begin{cor}\mlabel{cor:mpoaa}
Let $(A,A^*,-R_{\succ_A}^*,L_{\star_A}^*,
R_{\star_A}^*,-L_{\prec_A}^*,
-R_{\bullet_A}^*,-L_{\bullet_A}^*,-R_{\succ_{A^*}}^*,L_{\star_{A^*}}^*,
R_{\star_{A^*}}^*,-L_{\prec_{A^*}}^*,
-R_{\bullet_{A^*}}^*,-L_{\bullet_{A^*}}^*)$ be a matched pair of NS-algebras. Define binary operations $\star_A$ and $\star_{A^*}$ by Eq.~\meqref{eq:ast}, respectively. Then $(A,A^*,R_{\prec_A}^*,L_{\succ_A}^*,R_{\prec_A^*}^*,L_{\succ_A^*}^*)$ is a matched pair of associative algebras. 
\end{cor}
\begin{proof}
By Eq.~\meqref{eq:lrast}, we obtain
$$l_{\star_1}:=-R_{\succ_A}^*+R_{\star_A}^*-R_{\bullet_A}^*=R_{\prec_A}^*,\quad
r_{\star_1}:=L_{\star_A}^*-L_{\prec_A}^*-L_{\bullet_A}^*=L_{\succ_A}^*,$$ 
and
$$l_{\star_2}:=-R_{\succ_{A^*}}^*+R_{\star_{A^*}}^*-R_{\bullet_{A^*}}^*=R_{\prec_{A^*}}^*,\quad r_{\star_2}:=L_{\star_{A^*}}^*-L_{\prec_{A^*}}^*-L_{\bullet_{A^*}}^*=L^*_{\succ_{A^*}}.$$
Thus, by Proposition \mref{pro:mpoaa}, $(A,A^*,R_{\prec_A}^*,L_{\succ_A}^*,R_{\prec_{A^*}}^*,L_{\succ_{A^*}}^*)$ is a matched pair of associative algebras. 
\end{proof}

\begin{pro}\mlabel{pro:ac}
Let $(A,\delta_{\prec_A},\delta_{\succ_A},
\delta_{\bullet_A})$ be an NS-coalgebra. Then $(A,\delta_{\star})$ is a coassociative coalgebra, where $\delta_\star$ is defined by Eq.~\meqref{eq:deltastar}.
\end{pro}
\begin{proof}
By Definition \mref{def:ns-coalgebra}, we have 
\begin{eqnarray*}
	(\delta_{\star}\ot\id)\delta_{\star}
	&=&(\delta_{\prec}\ot\id)\delta_{\prec}+(\delta_{\succ}\ot\id)\delta_{\prec}+(\delta_{\bullet}\ot\id)\delta_{\prec}+(\delta_{\star}\ot\id)\delta_{\succ}+(\delta_{\star}\ot\id)\delta_{\bullet}\\
	&\stackrel{(\mref{eq:NSca04})}=&(\delta_{\prec}\ot\id)\delta_{\prec}+(\delta_{\succ}\ot\id)\delta_{\prec}+(\delta_{\star}\ot\id)\delta_{\succ}+(\id\ot\delta_{\bullet})\delta_{\succ}+(\id\ot\delta_{\star})\delta_{\bullet}\\
	&\stackrel{(\mref{eq:NSca02})}=&(\id\ot\delta_\star)\delta_\prec+(\delta_{\succ}\ot\id)\delta_{\prec}+(\delta_{\star}\ot\id)\delta_{\succ}+(\id\ot\delta_{\bullet})\delta_{\succ}+(\id\ot\delta_{\star})\delta_{\bullet}\\
	&\stackrel{(\mref{eq:NSca01})}=&(\id\ot\delta_\star)\delta_\prec+(\id\ot\delta_\prec)\delta_\succ+(\delta_{\star}\ot\id)\delta_{\succ}+(\id\ot\delta_{\bullet})\delta_{\succ}+(\id\ot\delta_{\star})\delta_{\bullet}\\
	&\stackrel{(\mref{eq:NSca03})}=&(\id\ot\delta_\star)\delta_\prec+(\id\ot\delta_\prec)\delta_\succ+(\id\ot\delta_\succ)\delta_\succ+(\id\ot\delta_{\bullet})\delta_{\succ}+(\id\ot\delta_{\star})\delta_{\bullet}\\
	&=&(\id\ot\delta_{\star})\delta_{\star}.
\end{eqnarray*}
Thus, $(A,\delta_{\star})$ is a coassociative coalgebra.
\end{proof}

\begin{defi}
	A {\bf dendriform D-bialgebra}~\cite{Bai1} is a quintuple
	$(A,\prec_A,\succ_A,\delta_{\prec_A},\delta_{\succ_A})$
	consisting of a dendriform algebra $(A,\prec_A,\succ_A)$
	and a dendriform coalgebra
	$(A,\delta_{\prec_A},\delta_{\succ_A})$
	such that, for all $x,y\in A$,
	\begin{align*}
		\delta_{\prec_A}(x\ast_A y)
		&=(\id\ot L_{\succ_A}(x))\delta_{\prec_A}(y)
		+(R_{\ast_A}(y)\ot\id)\delta_{\prec_A}(x),
		\\
		\delta_{\succ_A}(x\ast_A y)
		&=(R_{\prec_A}(y)\ot\id)\delta_{\succ_A}(x)
		+(\id\ot L_{\ast_A}(x))\delta_{\succ_A}(y),
		\\
		\delta_{\ast_A}(x\prec_A y)
		&=(\id\ot L_{\prec_A}(x))\delta_{\succ_A}(y)
		+(R_{\prec_A}(y)\ot\id)\delta_{\ast_A}(x),
	\\
		\delta_{\ast_A}(x\succ_A y)
		&=(R_{\succ_A}(y)\ot\id)\delta_{\prec_A}(x)
		+(\id\ot L_{\succ_A}(x))\delta_{\ast_A}(y),\\
		\big((\id\ot R_{\prec_A}(y))-(L_{\succ_A}(y)\ot\id)\big)\delta_{\ast_A}(x)
		&=\tau\big((L_{\prec_A}(x)\ot\id)\delta_{\prec_A}(y)
		-(\id\ot R_{\succ_A}(x))\delta_{\succ_A}(y)\big),
		\\
		\big((L_{\ast_A}(x)\ot\id)-(\id\ot R_{\prec_A}(x))\big)\delta_{\prec_A}(y)
		&=\tau\big((\id\ot R_{\ast_A}(y))\delta_{\succ_A}(x)
		-(L_{\succ_A}(y)\ot\id)\delta_{\succ_A}(x)\big).		
	\end{align*}
	Here $\ast_A,\ast_{A^*}$ are defined by
	Eq.~\meqref{eq:astoda},
	$\tau$ is the flip map, and
	$\delta_{\prec_A},\delta_{\succ_A},\delta_{\ast_A}$
	are the dual linear maps of
	$\prec_{A^*},\succ_{A^*},\ast_{A^*}$,
	respectively.
\end{defi}
By Remark~\mref{rmk:da} and
Definition~\mref{def:ns-bialgebra},
we obtain the following immediate consequence.
\begin{cor}\mlabel{cor:ippnslb}
	Let $(A,\prec_A,\succ_A,\bullet_A,
	\delta_{\prec_A},\delta_{\succ_A},
	\delta_{\bullet_A})$ be an NS-bialgebra
	such that $\bullet_A=0$ and
	$\delta_{\bullet_A}=0$.
	Then $(A,\prec_A,\succ_A,
	\delta_{\prec_A},\delta_{\succ_A})$
	is a dendriform D-bialgebra.
\end{cor}

\begin{lem}\mlabel{lem:aal}
	Let $(A,\prec,\succ,\bullet)$ be an NS-algebra.
	Then $(A\oplus A,\diamond)$ is an associative
	algebra, where $\diamond$ is defined by
	\begin{align}
		(x_1,x_2)\diamond(y_1,y_2)
		:=(x_1\star y_1,x_1\succ y_2),
		\quad (x_1,x_2),(y_1,y_2)\in A\oplus A.
		\mlabel{eq:aaol}
	\end{align}
\end{lem}

\begin{proof}
	For $(x_1,x_2)$, $(y_1,y_2)$
	and $(z_1,z_2)$ in $A\oplus A$,
	we have
	\begin{align*}
		((x_1,x_2)\diamond (y_1,y_2))\diamond (z_1,z_2)
		&=((x_1\star y_1)\star z_1,
		(x_1\star y_1)\succ z_2),\\
		(x_1,x_2)\diamond ((y_1,y_2)\diamond (z_1,z_2)
		&=(x_1\star(y_1\star z_1),
		x_1\succ(y_1\succ z_2)).
	\end{align*}
Since $\star$ is associative and
$(x\star y)\succ z=x\succ(y\succ z)$
by Eq.~\meqref{eq:nsa03},
the two expressions coincide, proving
the associativity of $\diamond$.
\end{proof}
\begin{rmk}
	An alternative associative algebra structure
	on $A\oplus A$ is given by
	\[
	(x_1,x_2)\diamond'(y_1,y_2):=(y_1\prec x_2,y_2\star x_2),\quad(x_1,x_2),(y_1,y_2)\in A\oplus A.
	\]
	Indeed, the associativity of $\diamond'$
	follows from the associativity of $\star$
	and Eq.~\meqref{eq:nsa01}.
\end{rmk}

\begin{defi}\cite{Bai1}
	Let $(A,\cdot)$ be an associative algebra. An antisymmetric bilinear form $\omega$ is called a {\bf Connes cocycle} on $(A,\cdot)$, if 
	$$\omega(xy,z)+\omega(yz,x)+\omega(zx,y)=0,\quad x,y,z\in A.$$
\end{defi}

\begin{lem}\mlabel{lem:nccoaal}
	Let $(A,\prec,\succ,\bullet,\frakB)$
	be a quadratic NS-algebra.
	Then $\omega$ is a nondegenerate Connes
	cocycle on $(A\oplus A,\diamond)$,
	where $\diamond$ is defined by
	Eq.~(\mref{eq:aaol}) and $\omega$ is given by
	\begin{align}
		\omega((x_1,x_2),(y_1,y_2))
		:=\frakB(x_1,y_2)+\frakB(x_2,y_1),
		\quad (x_1,x_2),(y_1,y_2)\in A\oplus A.
		\mlabel{eq:ndccol}
	\end{align}
\end{lem}

\begin{proof}
	By Lemma~\mref{lem:aal},
	$(A\oplus A,\diamond)$ is an associative
	algebra. Since $\frakB$ is nondegenerate
	and antisymmetric, so is $\omega$.
	First, by Eqs.~(\mref{eq:ibf2})--(\mref{eq:ibf3})
	and the antisymmetry of $\frakB$, we obtain
	\begin{align}
		\frakB(x\star y,z)
		&=\frakB(x,y\star z)
		+\frakB(y,z\star x)
		-\frakB(x,y\bullet z)\notag\\
		&=\frakB(x,y\prec z)
		+\frakB(x,y\succ z)
		+\frakB(y,z\star x)\notag\\
		&=\frakB(x,y\succ z). \mlabel{eq:xsyz}
	\end{align}
For $(x_1,x_2)$, $(y_1,y_2)$
	and $(z_1,z_2)$ in $A\oplus A$,
	we have
	\begin{align*}
		&\omega((x_1,x_2)\diamond (y_1,y_2),(z_1,z_2))
		+\omega((y_1,y_2)\diamond (z_1,z_2),(x_1,x_2))
		+\omega((z_1,z_2)\diamond (x_1,x_2),(y_1,y_2))\\
		={}&\frakB(x_1\star y_1,z_2)
		+\frakB(x_1\succ y_2,z_1)
		+\frakB(y_1\star z_1,x_2)
		+\frakB(y_1\succ z_2,x_1)\\
		&+\frakB(z_1\star x_1,y_2)
		+\frakB(z_1\succ x_2,y_1)\\
		\stackrel{\meqref{eq:xsyz}}{=}&\frakB(x_1,y_1\succ z_2)
		-\frakB(z_1,x_1\succ y_2)
		+\frakB(y_1,z_1\succ x_2)
		-\frakB(x_1,y_1\succ z_2)\\
		&+\frakB(z_1,x_1\succ y_2)
		-\frakB(y_1,z_1\succ x_2)\\
		={}&0.
	\end{align*}
	Hence $\omega$ is a nondegenerate
	Connes cocycle on $(A\oplus A,\diamond)$.
\end{proof}

\begin{defi}\cite{Bai1}
Let $(A,\cdot)$ be an associative algebra, and let $\frakB_d$ be an antisymmetric bilinear form given by Eq.~\meqref{eq:Bd1}. $(A,\frakB_d)$ is called a {\bf double construction of the Connes cocycle}, if it satisfies the conditions
\begin{enumerate}
	\item\mlabel{it:dcotcc1}
	$A=A_1\oplus A^*_1$ as the direct sum of vector spaces;
	
	\item\mlabel{it:dcotcc2}
	$A_1$ and $A^*_1$ are associative subalgebras
	of $A$;
	
	\item\mlabel{it:dcotcc3}
	$\frakB_d$  is a  Connes cocycle on $(A,\cdot)$.
\end{enumerate}
\end{defi}

\begin{pro}\mlabel{pro:dcotcc}
	Let $\bigl((D:=A\oplus A^*,
	\prec_{\bowtie},\succ_{\bowtie},
	\bullet_{\bowtie},\frakB_d),A,A^*\bigr)$
	be a Manin triple of NS-algebras.
	Then $(D\oplus D,\diamond_d,\omega_d)$
	is a double construction of the Connes cocycle,
	where $\diamond_d$ and $\omega_d$
	are defined by Eqs.~(\mref{eq:aaol})
	and (\mref{eq:ndccol}), respectively,
	with $A$ replaced by $D$.
\end{pro}

\begin{proof}
	Set $A_1:=A\oplus A$ and
	$B_1:=A^*\oplus A^*$, viewed as
	subspaces of $D\oplus D$ via the
	natural embeddings.
	Then $D\oplus D=A_1\oplus B_1$
	as vector spaces.
	Since $A$ and $A^*$ are NS-subalgebras
	of $D$, Lemma~\mref{lem:aal} implies
	that $A_1$ and $B_1$ are associative
	subalgebras of $(D\oplus D,\diamond_d)$.
	
	Since $(D,\prec_{\bowtie},\succ_{\bowtie},
	\bullet_{\bowtie},\frakB_d)$ is a
	quadratic NS-algebra,
	Lemma~\mref{lem:nccoaal} implies
	that $\omega_d$ is a nondegenerate
	Connes cocycle on
	$(D\oplus D,\diamond_d)$.
	Define a linear isomorphism
	$J:B_1\to A_1^*$ by
	\[
	\langle J(\xi,\eta),(x,y)\rangle
	=\langle\eta,x\rangle+\langle\xi,y\rangle.
	\]
	For $a_i=(x_i,y_i)\in A_1\subseteq D\oplus D$ and
$b_i=(\xi_i,\eta_i)\in B_1\subseteq D\oplus D$, we have
\[
a_i+b_i=(x_i+\xi_i,y_i+\eta_i)
\quad\text{in }D\oplus D.
\]
Hence, by the definition of $\omega_d$,
\begin{align*}
	\omega_d(a_1+b_1,a_2+b_2)
	&=\frakB_d(x_1+\xi_1,y_2+\eta_2)
	+\frakB_d(y_1+\eta_1,x_2+\xi_2)\\
	&=-\langle x_1,\eta_2\rangle
	-\langle y_1,\xi_2\rangle+\langle\xi_1,y_2\rangle
	+\langle\eta_1,x_2\rangle\\
	&=-\langle a_1,J(b_2)\rangle
	+\langle J(b_1),a_2\rangle.
\end{align*}
	Thus, under the identification
	$B_1\cong A_1^*$ induced by $J$,
	$\omega_d$ coincides with the
	canonical antisymmetric bilinear form
	on $A_1\oplus A_1^*$.
	Consequently,
	$(D\oplus D,\diamond_d,\omega_d)$
	is a double construction of the
	Connes cocycle.
\end{proof}

\begin{thm}\mlabel{thm:ns-to-dendriform}
	Every NS-bialgebra on a vector space $A$ induces
	a dendriform D-bialgebra structure
	on $A\oplus A$.
\end{thm}

\begin{proof}
	Let $(A,\prec_A,\succ_A,\bullet_A ,
	\delta_{\prec_A},\delta_{\succ_A},
	\delta_{\bullet_A} )$
	be an NS-bialgebra. By Theorem \mref{thm:equiv-ns}, $\bigl((D:=A\oplus  A^*,\prec_{\bowtie},\succ_{\bowtie},\bullet_{\bowtie},\frakB_d),A,A^*\bigr)$
	is a Manin triple of NS-algebras. Thus, by Proposition \mref{pro:dcotcc}, $(D\oplus D=A_1\oplus A^*_1,\diamond_d,\omega_d)$ is a double construction of the Connes cocycle, where $A_1:=A\oplus A$, $A^*_1:=A^*\oplus A^*$, $\diamond_d$ and $\omega_d$ are given by Eq.~(\mref{eq:aaol}) and Eq.~(\mref{eq:ndccol}), respectively. Therefore, by \cite[Theorem 4.2.4]{Bai1}, there is a dendriform D-bialgebra structure on $A_1$.
\end{proof}

	\noindent {\bf Acknowledgments}: This work was supported by the National Natural Science Foundation of China (12461002), and the Jiangxi Provincial Natural Science Foundation (20224BAB201003).
	The first author thanks the Chern Institute of Mathematics at Nankai University for hospitality.
	
	\noindent
	{\bf Declaration of interests. } The authors have no conflicts of interest to disclose.
	
	\noindent
	{\bf Data availability. } Data sharing is not applicable as no data were created or analyzed.\\

\end{document}